\documentclass{amsart}

\usepackage{amsmath}
\usepackage{amssymb}
\usepackage{amsthm}

\newtheorem{thm}{Theorem}[section]
\newtheorem{prop}[thm]{Proposition}
\newtheorem{lem}[thm]{Lemma}

\theoremstyle{definition}
\newtheorem{dfn}[thm]{Definition}

\theoremstyle{remark}
\newtheorem{rem}[thm]{Remark}
\newtheorem*{acknowledgement}{Acknowledgements}
\newtheorem*{aiuse}{AI use statement}

\newcommand{\R}{\mathbb{R}}

\newcommand{\Aff}{\mathrm{Aff}}

\newcommand{\Mag}{\mathrm{Mag}}
\newcommand{\1}{1\!\!1}
\newcommand{\SSL}{\mathrm{SSL}}

\title[Excess inertia of finite metric spaces]
{Excess inertia of finite metric spaces}

\author{Kiyonori Gomi}

\address{
Department of Mathematics, 
Institute of Science Tokyo,
2-12-1 Ookayama, Meguro-ku, Tokyo, 152-8551, Japan.}

\email{kgomi@math.titech.ac.jp}



\begin{document}

\begin{abstract}
This paper classifies the possible inertias that are realized by an excess matrix of a finite metric space with at most nine points. For the classification, we make use of relationships between finite metric spaces and non-obtuse and acute sets in classical geometry. As an application, we prove that the small-scale limit of the magnitude function of any six-point metric space is at least one whenever this limit exists.
\end{abstract}

\maketitle

\tableofcontents

\section{Introduction}
\label{sec:introduction}

Associated to a finite metric space are the so-called \textit{excess matrices}. To describe them, let $d$ be a metric on an $n$-point set $X = \{ 0, 1, \ldots, n-1 \}$. Then, for each $k \in X$, the $k$th excess matrix $E^{(k)}$ is the $(n-1) \times (n-1)$ symmetric matrix defined by $E^{(k)} = ( d_{k, i} + d_{k, j} - d_{i, j} )_{i, j \in X \backslash \{ k \}}$. It is easy to see that the excess matrices of $X = (X, d)$ have the same inertia, which will be called the \textit{excess inertia} of $X$ in this paper. We write $(n_+, n_-, n_0)$ for the inertia, where $n_+, n_-$ and $n_0$ are respectively the numbers of positive, negative and zero eigenvalues of $E^{(k)}$ counted with multiplicity.

A natural problem would be to classify the inertias $(n_+, n_-, n_0)$ that occur as the excess inertia of an $n$-point metric space. The excess inertia of a finite metric space $(X, d)$ should be distinguished from the inertia of the distance matrix defined by $D = ( d_{i, j} )_{i, j \in X}$. The classification problem of excess inertia may appear to be classical. However, we are not aware of a systematic study of this problem. 

It is well-known that a metric space is of \textit{negative type} \cite{DL} if and only if its excess matrix is positive semidefinite. Also, it is of \textit{strictly negative type} \cite{HLMT} if and only if the excess matrix is positive definite. Furthermore, any finite metric space consisting of less than or equal to four points is of negative type, and is of strictly negative type unless it is a ``metric rectangle'' (see \S\S\ref{subsec:schoenberg_embedding} for its definition) \cite{HLMT}, which has the excess inertia $(2, 0, 1)$.  These results lead one to classify the excess inertias of metric spaces consisting of $n \le 4$ points. For example, the possible excess inertias of $4$-point metric spaces are precisely as follows:
$$
\begin{array}{|c|c|c|}
\hline
\mbox{$4$-point} & n_0 = 0 & n_0 = 1 \\
\hline
n_- = 0 & (3, 0, 0) & (2, 0, 1) \\
\hline
\end{array}
$$
A $5$-point metric space need not be of negative type. An unpublished work of Asao and Taira classifies the possible excess inertias of $5$-point metric spaces:
$$
\begin{array}{|c|c|c|}
\hline
\mbox{$5$-point} & n_0 = 0 & n_0 = 1 \\
\hline
n_- = 0 & (4, 0, 0) & (3, 0, 1) \\
\hline
n_- = 1 & (3, 1, 0) & \\
\hline
\end{array}
$$

The main result of this paper is the classification of excess inertias of metric spaces consisting of $n \le 9$ points. The following theorem summarizes the results from Theorem \ref{thm:6_point}, Theorem \ref{thm:7_point}, Theorem \ref{thm:8_point} and Theorem \ref{thm:9_point} in the body of this paper.

\begin{thm} \label{ithm:excess_inertia}
For $n = 6, 7, 8, 9$, the possible excess inertias of $n$-point metric spaces are precisely as follows:
$$
\begin{array}{|c|c|c|c|}
\hline
\mbox{$6$-point} & n_0 = 0 & n_0 = 1 & n_0 = 2 \\
\hline
n_- = 0 & (5, 0, 0) & (4, 0, 1) & (3, 0, 2) \\
\hline
n_- = 1 & (4, 1, 0) & & \\
\hline
\end{array}
$$
$$
\begin{array}{|c|c|c|c|c|}
\hline
\mbox{$7$-point} & n_0 = 0 & n_0 = 1 & n_0 = 2 & n_0 = 3 \\
\hline
n_- = 0 & (6, 0, 0) & (5, 0, 1) & (4, 0, 2) & (3, 0, 3)\\
\hline
n_- = 1 & (5, 1, 0) & (4, 1, 1) & &  \\
\hline
n_- = 2 & (4, 2, 0) & & & \\
\hline
\end{array}
$$
$$
\begin{array}{|c|c|c|c|c|c|}
\hline
\mbox{$8$-point} & n_0 = 0 & n_0 = 1 & n_0 = 2 & n_0 = 3 & n_0 = 4 \\
\hline
n_- = 0 & (7, 0, 0) & (6, 0, 1) & (5, 0, 2) & (4, 0, 3) & (3, 0, 4) \\
\hline
n_- = 1 & (6, 1, 0) & (5, 1, 1) & (4, 1, 2) & & \\
\hline
n_- = 2 & (5, 2, 0) & (4, 2, 1) & & & \\
\hline
n_- = 3 & (4, 3, 0) & & & & \\
\hline
\end{array}
$$
$$
\begin{array}{|c|c|c|c|c|c|c|c|}
\hline
\mbox{$9$-point} & n_0 = 0 & n_0 = 1 & n_0 = 2 & n_0 = 3 & n_0 = 4  \\
\hline
n_- = 0 & (8, 0, 0) & (7, 0, 1) & (6, 0, 2) & (5, 0, 3) & (4, 0, 4) \\
\hline
n_- = 1 & (7, 1, 0) & (6, 1, 1) & (5, 1, 2) & (4, 1, 3) &  \\
\hline
n_- = 2 & (6, 2, 0) & (5, 2, 1) & (4, 2, 2) & & \\
\hline
n_- = 3 & (5, 3, 0) & (4, 3, 1) &  & & \\
\hline
n_- = 4 & (4, 4, 0)&  & & & \\
\hline
\end{array}
$$
\end{thm}

The proof of the above result relies on relationships between the classification of excess inertia and non-obtuse and acute sets in classical geometry. Recall that a set $S$ in the Euclidean space $\R^d$ of dimension $d$ is called a \textit{non-obtuse set} (resp.\ \textit{an acute set}) if the angle $\angle u v w$ is non-obtuse (resp.\ acute) for any three distinct points $u, v, w \in S$. In terms of the standard inner product $\langle \ \ , \  \ \rangle$ on $\R^d$, the non-obtuse and acute conditions are expressed as $\langle u - v, w - v \rangle \ge 0$ and $\langle u - v, w - v \rangle > 0$, respectively. For non-negative integers $p, q$, let $\R^{p, q}$ denote the $(p+q)$-dimensional real vector space equipped with the inner product
$$
\langle (x_1, \ldots, x_{p+q}), (y_1, \ldots, y_{p+q}) \rangle_{p, q} 
= \sum_{i = 1}^p x_iy_i - \sum_{i = p+1}^{p+q} x_iy_i.
$$
In view of the above expressions of the non-obtuse and acute conditions in $\R^d$, we call a set $S \subset \R^{p, q}$ non-obtuse (resp.\ acute) if $\langle u - v, w - v \rangle_{p, q} \ge 0$ (resp.\ $\langle u - v, w - v \rangle_{p, q} > 0$) for any three distinct points $u, v, w \in S$. Then, thanks to a classical idea of Schoenberg \cite{Sch}, if a metric $d$ on a finite set $X = \{ 0, 1, \ldots, n - 1 \}$ consisting of $n \ge 3$ points has the excess inertia $(p, q, n - p - q - 1)$, then we can express the metric as $d_{i, j} = \lVert v_i - v_j \rVert^2_{p,q}$ by means of a non-obtuse set $\{ v_0, \ldots, v_{n-1} \}$ in $\R^{p, q}$ (Proposition \ref{prop:metric_and_non_obtuse_set}). Then the following relationships will be used in this paper: 

\begin{itemize}
\item
By the result of Danzer and Gr\"{u}nbaum \cite{DG} solving a problem of Erd\H{o}s \cite{E}, any non-obtuse set in $\R^p = \R^{p, 0}$ contains at most $2^p$ points. As a consequence, if $(p, 0, n - p - 1)$ is the excess inertia of an $n$-point metric space, then $p_- \le p \le n - 1$, where $p_- = \min \{ p |\ n \le 2^p \} = \lceil \log_2 n \rceil$. For example, $(3, 0, 5)$ cannot be the excess inertia of a $9$-point metric space. 

\item
If the affine span of an $n$-point acute set $S$ in $\R^p = \R^{p, 0}$ is the whole space $\R^p$, then we get an $n$-point acute set in $\R^{p, q}$ for any $q = 0, 1, \ldots, n - p - 1$ (Lemma \ref{lem:acute}). It is known by Gerencs\'{e}r and Harangi \cite{GH} that there exists a $(2^{p-1} + 1)$-point acute set in $\R^p$. Thus, for example, there exists a $9$-point metric space whose excess inertia is $(4, 4, 0)$.

\item
In order to establish the main theorem, the most non-trivial part of the proof is the non-existence of a $6$-point metric space of excess inertia $(3, 1, 1)$ (Proposition \ref{prop:6_point_(3,1,1)}). Our proof of this non-existence relies on the classification of $6$-point non-obtuse sets in $\R^3$ obtained by Croft \cite{Croft}.

\end{itemize}

For the time being, our result about the excess inertia of a $10$-point metric space is partial (Proposition \ref{prop:10_point}). For its completion, the classification of $10$-point non-obtuse sets in $\R^4$ may play the central role, in view of our argument to classify the excess inertias of $6$-point metric spaces.

\medskip

We have two applications of Theorem \ref{ithm:excess_inertia}.

\smallskip

The first application concerns a generalization of the classical non-obtuse set problem in Euclidean spaces to $\R^{p, q}$ for small $p, q$. Let us consider a set $S$ in $\R^{p, q}$ satisfying the following three conditions.

\begin{enumerate}
\item[(i)]
We have $\lVert u - v \rVert^2_{p, q} > 0$ for distinct $u, v \in S$.

\item[(ii)]
We have $\langle u - v, w - v \rangle_{p, q} \ge 0$ for distinct $u, v, w \in S$.

\item[(iii)]
The affine span of $S$ is the whole space: $\Aff(S) = \R^{p, q}$.

\end{enumerate}
It turns out that $S$ consists of a finite number of points. Indeed, the idea of the proof of Danzer-Gr\"{u}nbaum works for non-obtuse sets in $\R^{p, q}$, and $S$ contains at most $2^{p+q}$ points (Lemma \ref{lem:finiteness}). However, $2^{p+q}$ cannot be generally attained in the case where $q > 0$. Let $g(p, q)$ be the maximal number of points contained in such an $S \subset \R^{p, q}$ as above. We set $g(p, q) = 0$ if no such set exists. By the correspondence of a metric space of excess inertia $(n_+, n_-, n_0)$ and a non-obtuse set in $\R^{n_+, n_-}$ whose affine span is the whole space, we get the following result.

\begin{thm} \label{ithm:non_obtuse_problem}
For small $p, q$, the numbers $g(p, q)$ are as follows:
$$
\begin{array}{|c|c|c|c|c|c|c|}
\hline
g(p, q) & q = 0 & q = 1 & q = 2 & q \ge 3 \\
\hline
p = 0 & 1 & 0 & 0 & 0 \\
\hline
p = 1 & 2 & 0 & 0 & 0 \\
\hline
p = 2 & 4 & 0 & 0 & 0 \\
\hline
p = 3 & 8 & 5 & 0 & 0 \\
\hline
\end{array}
$$
\end{thm}

The theorem above follows from Theorem \ref{thm:maximal_points_in_non_obtuse_set_in_Rpq} in the body of this paper. Note that $g(p, 0) = 2^p$ is nothing but the theorem of Danzer and Gr\"{u}nbaum \cite{DG}. 

\smallskip

The second application concerns Leinster's notion of \textit{magnitude} \cite{L}. The magnitude $\Mag(X)$ of a finite metric space $X = (X, d)$ is an isometric invariant which can be viewed as counting the points in $X$ while taking the metric $d$ into account. For a positive real number $t$, we let $tX = (X, td)$ be the metric space scaled by $t$. The limit of $\Mag(tX)$ as $t \to 0$ is, if it exists, called the \textit{small-scale limit} of the magnitude (function) of $X$, and we put $\SSL(X) = \lim_{t \to 0}\Mag(tX)$. The small-scale limit is systematically studied in a work of Roff and Yoshinaga \cite{RY}, and a use of some classical results about the magnitude shows that $\SSL(X) = 1$ for any non-empty metric space consisting of $n \le 4$ points. An unpublished result of Asao and Taira states that $\SSL(X) \ge 1$ for any $5$-point metric space $X$ which admits a small-scale limit. Essentially, their proof makes use of the classification of excess inertias of $5$-point metric spaces. Generalizing their idea based on the classification of excess inertias of $6$-point metric spaces in this paper, we get the following result (which will be proved as Theorem \ref{thm:6_point_SSL}).

\begin{thm} \label{ithm:SSL_6_point}
We have $\SSL(X) \ge 1$ for any $6$-point metric space $X$ which admits a small-scale limit.
\end{thm}

The classification of excess inertia was originally motivated by this application.

\bigskip

This paper is organized as follows. In \S\ref{sec:preliminaries}, we recall some basic facts about metric spaces and non-obtuse sets and acute sets in Euclidean spaces. Some constructions to be used are also summarized in this section. In \S\ref{sec:classification_of_excess_inertia}, we study the excess inertia of metric spaces consisting of $n \le 10$ points. In \S\ref{sec:application}, our classification of excess inertia is applied to the non-obtuse set problem in $\R^{p, q}$ and the small-scale limit of magnitude.

\begin{acknowledgement}
The author thanks Yasuhiko Asao and Kouichi Taira for sharing their results and helpful discussions.
\end{acknowledgement}

\begin{aiuse}
The author used GPT-6 Pro to assist in developing the proof of Proposition \ref{prop:6_point_(3,1,1)}. The author has verified the final argument and computations and takes full responsibility for the mathematical content of this paper.
\end{aiuse}

\section{Preliminaries}
\label{sec:preliminaries}

\subsection{Excess inertia}

Let $X = \{ 0, 1, \ldots, n-1 \}$ be a finite set consisting of $n$ points. A map $d : X \times X \to [0, \infty)$ is said to be a \textit{semimetric} when (i) and (ii) below are satisfied, and is said to be a \textit{metric} when (i), (ii) and (iii) are satisfied.
\begin{enumerate}
\item[(i)]
$d_{x, y} = 0$ if and only if $x = y \in X$. 

\item[(ii)]
$d_{x, y} = d_{y, x}$ for all $x, y \in X$.

\item[(iii)]
$d_{x, y} + d_{y, z} \ge d_{x, z}$ for all $x, y, z \in X$.

\end{enumerate}
As is well known, the inequality in (iii) is called the triangle inequality. When $d$ is a (semi)metric on $X$, we call $d$ the distance (function). The pair $(X, d)$ is called a (semi)metric space. We often write $X$ for $(X, d)$.

The \textit{distance matrix} $D \in M(n, \R)$ of $(X, d)$ (or $d$) is the symmetric matrix defined by $D = (d_{i, j})_{i, j \in X}$, where $M(n, \R)$ is the set of $n \times n$ matrices with real entries. Suppose $n > 0$. For $k \in X$, the $k$th \textit{excess matrix} $E^{(k)} \in M(n-1, \R)$ of $(X, d)$ (or $d$) is defined by $E^{(k)} = ( d_{i, k} + d_{j, k} - d_{i, j} )_{i, j \in X \backslash \{ k \}}$. For any $k$, the inertia of $E^{(k)}$ is identified with that of the quadratic form $- \langle v, D v \rangle$ for $v \in \Pi_0 = \{ v = (v_i)_{i \in X} \in \R^n |\ \langle \1, v \rangle = 0 \}$, where $\langle \ , \ \rangle$ is the standard inner product and $\1$ the vector whose entries are $1$. Hence all the excess matrices $E^{(k)}$ have the same inertia, and the following notion is well-defined.

\begin{dfn}
Let $(X, d)$ be a non-empty finite semimetric space. We define the \textit{excess inertia} of $(X, d)$ or $d$ as the inertia of an excess matrix of $(X, d)$.
\end{dfn}

We will write $(n_+, n_-, n_0)$ for the inertia of an excess matrix $E \in M(n-1, \R)$, where $n_+$, $n_-$ and $n_0$ are respectively the numbers of positive, negative and zero eigenvalues of $E$ counted with multiplicity, and hence $n_+ + n_- + n_0 = n - 1$. Because the diagonal entries of $E$ are positive, $n_+ > 0$ if $n \ge 2$.

\medskip

\begin{rem} \label{rem:negative_type}
In metric geometry, the notion of \textit{negative type} \cite{DL} is classical. For $n \ge 2$, one can easily see that a metric $d$ on an $n$-point set $X$ is of negative type if and only if an excess matrix is positive semidefinite, and hence the excess inertia of $d$ is $(p, 0, n - p - 1)$ for some $p \in \{ 1, 2, \ldots, n - 1 \}$. It is also easy to see that $d$ is \textit{strictly negative type} \cite{HLMT} if and only if an excess matrix is positive definite, and hence the excess inertia of $d$ is $(n-1, 0, 0)$.
\end{rem}

\medskip

The following lemma is a consequence of the Cauchy interlace theorem \cite{HJ}.

\begin{lem} \label{lem:possible_inertia}
Let $(n_+, n_-, n_0)$ be the excess inertia of an $n$-point metric space $X$. Then the excess inertia of any $(n+1)$-point metric space containing $X$ as a subspace must be one of the following
\begin{align*}
&(n_+ + 1, n_-, n_0), &
&(n_+, n_- + 1, n_0), &
&(n_+, n_-, n_0 + 1), &
&(n_+ + 1, n_- + 1, n_0 - 1), 
\end{align*}
where the last case occurs only when $n_0 \ge 1$.
\end{lem}

\begin{proof}
Let $X_+ = \{ 0, 1, \ldots, n \}$ be an $(n+1)$-point metric space which contains $X = \{ 0, 1, \ldots, n-1 \}$. Let $E$ and $E_+$ be the $0$th excess matrices of $X$ and $X_+$, respectively. Then $E \in M(n-1, \R)$ is a principal submatrix of $E_+ \in M(n, \R)$. Let $\lambda_1 \le \cdots \le \lambda_{n-1}$ denote the eigenvalues of $E$ counted with multiplicity. Similarly, $\lambda_1^+ \le \cdots \le \lambda^+_{n}$ denote the eigenvalues of $E_+$. Then the Cauchy interlace theorem \cite{HJ} implies that 
$$
\lambda_1^+ \le \lambda_1 \le
\lambda_2^+ \le \lambda_2 \le
\lambda_3^+ \le
\cdots
\le
\lambda_{n-1}^+ \le \lambda_{n-1} \le \lambda_n^+.
$$
Now, using this relationship of eigenvalues, we consider the possibilities of the inertia of $E_+$, provided that the inertia of $E$ is $(n_+, n_-, n_0)$. 
\begin{itemize}
\item
In the case where $n_0 = 0$, the matrix $E_+$ has at least $n_+$ positive eigenvalues and at least $n_-$ negative eigenvalues. There is one more eigenvalue of $E_+$, which is either positive, negative or zero. 

\item
In the case where $n_0 \ge 1$, the matrix $E_+$ has at least $n_+$ positive eigenvalues, at least $n_-$ negative eigenvalues and at least $n_0 - 1$ zero eigenvalues. There are two more eigenvalues of $E_+$, whose combinations of signs are either $(+, 0)$, $(-, 0)$, $(0, 0)$ or $(+, -)$.

\end{itemize}
Summarizing this consideration, we get the possible inertia of $E_+$ as stated.
\end{proof}

\subsection{Schoenberg embedding}
\label{subsec:schoenberg_embedding}

For our analysis of the excess inertia, a particular realization of a semimetric, which we call the ``Schoenberg embedding'', will play the key role. To explain this notion, suppose that we are given non-negative integers $p, q$. For $v = (v_i), w = (w_i) \in \R^{p+q}$, we put
$$
\langle v, w \rangle_{p, q}
= \sum_{i = 1}^p v_iw_i - \sum_{i = p+1}^{p+q} v_iw_i.
$$
If we use a diagonal matrix $I_{p, q} = \mathrm{diag}(\overbrace{1, \ldots, 1}^p, \overbrace{-1, \ldots, -1}^q)$ of inertia $(p, q, 0)$, then $\langle v, w \rangle_{p, q} = \langle v, I_{p, q} w \rangle$ by using the standard inner product on $\R^{p + q}$. Note that $\lVert v \rVert^2_{p, q} = \langle v, v \rangle_{p, q}$ can be negative if $pq \neq 0$. We will write $\R^{p, q}$ to mean the vector space $\R^{p+q}$ equipped with $\langle \ , \ \rangle_{p, q}$.

Given a set of vectors $S = \{ v_0, \ldots, v_m \} \subset \R^{p, q}$, we define its \textit{affine span} by 
$$
\Aff(S) = 
\bigg\{
\sum_{i = 0}^m \alpha_i v_i \in \R^{p, q} \bigg|\
\alpha_0, \ldots, \alpha_m \in \R, \
\sum_{i = 0}^m \alpha_i = 1
\bigg\}.
$$

\begin{prop}[\cite{Sch}] \label{prop:Schoenberg_embedding}
Let $d$ be a semimetric on a finite set $X = \{ 0, 1, \ldots, n - 1 \}$ consisting of $n > 1$ points. Suppose that the inertia of its $0$th excess matrix $E^{(0)} \in M(n-1, \R)$ is $(p, q, n - p - q - 1)$. Then there exists a set of vectors $S = \{ v_0, v_1, \ldots, v_{n-1} \} \subset \R^{p, q}$ satisfying the following properties.
\begin{itemize}
\item[(a)]
$\Aff(S) = \R^{p, q}$.

\item[(b)]
We have $d_{i, j} = \lVert v_i - v_j \rVert^2_{p, q}$ for any $i, j \in X$.

\item[(c)]
$E^{(0)}/2$ agrees with the Gram matrix of the vectors $v_1 - v_0, \ldots, v_{n-1} - v_0$ in $\R^{p, q}$
$$
\frac{1}{2}E^{(0)} 
= \mathrm{Gram}(v_1 - v_0, \ldots, v_{n-1} - v_0) 
= ( \langle v_i - v_0, v_j - v_0 \rangle_{p, q} )_{i, j = 1}^{n-1}.
$$

\end{itemize}
\end{prop}

\begin{proof}
The proof is essentially given by Schoenberg \cite{Sch}: Let $\lambda_1, \ldots, \lambda_p > 0$ be the positive eigenvalues of $E^{(0)}/2$, and $- \lambda_{p+1}, \ldots, - \lambda_{p+q} < 0$ the negative eigenvalues of $E^{(0)}/2$. Let $\lambda_{p+q+1} = \cdots = \lambda_{n-1} = 0$ denote the remaining zero eigenvalues of $E^{(0)}/2$. Let $w_1, \ldots, w_{n-1} \in \R^{n-1}$ be the corresponding eigenvectors, which are orthonormalized in the standard inner product: $\langle w_i, w_j \rangle = \delta_{i, j}$. We here introduce $W, L, J \in M(n-1, \R)$ by
\begin{align*}
W &= {}^t(w_1, \ldots, w_{n-1}), &
L &=
\mathrm{diag}(
\sqrt{\lambda_1}, \ldots, \sqrt{\lambda_p},
\sqrt{\lambda_{p+1}}, \ldots, \sqrt{\lambda_{p+q}}, 
0, \ldots, 0),
\\
& &
J &= 
\mathrm{diag}(
\underbrace{1, \ldots, 1}_p,
\underbrace{-1, \ldots, -1}_q,
0, \ldots, 0
).
\end{align*}
We then have
\begin{align*}
\frac{1}{2}E^{(0)}
&= {}^tW 
\mathrm{diag}(
\lambda_1, \ldots, \lambda_p,
-\lambda_{p+1}, \ldots, -\lambda_{p+q}, 
0, \ldots, 0
)
W
\\
&=
{}^tW LJL W
=
{}^t(LW) J (LW),
\end{align*}
where ${}^t$ stands for the transpose. In view of the above expression, we define $w'_1, \ldots, w'_{n-1} \in \R^{n-1}$ by $(w'_1, \ldots, w'_{n-1}) = LW$. Using the projection $\pi : \R^{n-1} \to \R^{p + q}$ omitting the last $n-p-q-1$ entries, we further define $v_1, \ldots, v_{n-1} \in \R^{p+q}$ by $v_i = \pi(w'_i)$. In general, $\langle \pi(u), \pi(u') \rangle_{p, q} = \langle u, J u' \rangle$ for $u, u' \in \R^{n-1}$. Thus, 
\begin{align*}
\frac{1}{2}E^{(0)}
&=
( \langle w'_i, J w'_j \rangle )_{i, j = 1}^{n-1}
=
( \langle v_i, v_j \rangle_{p, q} )_{i, j = 1}^{n-1}
=
\mathrm{Gram}(v_1, \ldots, v_{n-1}).
\end{align*}
Now, the set of vectors $S = \{ v_0, v_1, \ldots, v_{n-1} \} \subset \R^{p, q}$ is given by adding $v_0 = 0$. Then (c) is clearly satisfied. It is easy to see that (c) is equivalent to (b). By construction, the vectors $v_1, \ldots, v_{n-1}$ span $\R^{p+q} = \R^{p, q}$, so that (a) holds.
\end{proof}

Generalizing the classical notions of non-obtuse sets and acute sets in a Euclidean space \cite{E}, we make the following definition.

\begin{dfn}
Let $S = \{ v_i \}_{i \in I} \subset \R^{p, q}$ be a set of (distinct) vectors. We say $S$ is  \textit{non-obtuse} if $\langle v_i - v_j, v_k - v_j \rangle_{p, q} \ge 0$ for all  pairwise distinct $i, j, k \in I$. We also say $S$ is \textit{acute} if $\langle v_i - v_j, v_k - v_j \rangle_{p, q} > 0$ for all pairwise distinct $i, j, k \in I$.
\end{dfn}

\begin{prop} \label{prop:metric_and_non_obtuse_set}
Let $n \ge 3$, and let $X = \{ 0, 1, \ldots, n-1 \}$ be an $n$-point set. For any metric $d$ on $X$ of excess inertia $(p, q, n - p - q - 1)$, there exists a non-obtuse set $S = \{ v_0, \ldots, v_{n-1} \}$ in $\R^{p, q}$ such that $\Aff(S) = \R^{p, q}$ and $d_{i, j} = \lVert v_i - v_j \rVert^2_{p, q}$ for all $i, j \in X$. Conversely, any non-obtuse set $S = \{ v_0, \ldots, v_{n-1} \}$ such that $\Aff(S) = \R^{p, q}$ defines a metric $d$ on $X$ of excess inertia $(p, q, n - p - q - 1)$ by $d_{i, j} =  \lVert v_i - v_j \rVert^2_{p, q}$.
\end{prop}

\begin{proof}
Given a metric $d$ on $X$ of excess inertia $(p, q, n - p - q - 1)$, Proposition \ref{prop:Schoenberg_embedding} provides us with a set $S = \{ v_0, \ldots, v_{n-1} \} \subset \R^{p, q}$ such that $\Aff(S) = \R^{p, q}$ and  $d_{i, j} = \lVert v_i - v_j \rVert^2_{p, q}$ for all $i, j \in X$. For any $i, j, k \in X$, we have
$$
d_{i, j} + d_{j, k} - d_{i, k}
= 2 \langle v_i - v_j, v_k - v_j \rangle_{p, q} \ge 0,
$$
because $d$ is a metric. Hence $S$ is non-obtuse. Conversely, suppose that $S = \{ v_0, \ldots, v_{n-1} \}$ is a non-obtuse set in $\R^{p, q}$, and put $d_{i, j} = \lVert v_i - v_j \rVert^2_{p, q}$ for $i, j \in X$. For a distinct pair $i, j \in X$, there exists $k \in X$ such that $k \neq i$ and $k \neq j$. We then have
\begin{align*}
\lVert v_i - v_j \rVert^2_{p, q}
&=
\langle (v_i - v_k) + (v_k - v_j), v_i - v_j \rangle_{p, q}
\\
&=
\langle v_i - v_k, v_i - v_j \rangle_{p, q}
+
\langle v_k - v_j, v_i - v_j \rangle_{p, q}
\ge 0.
\end{align*}
Assume here that $d_{i, j} = 0$. Then, for any $k$ distinct from both $i$ and $j$, the estimate above implies $\langle v_i - v_k, v_i - v_j \rangle_{p, q} = 0$. It follows that $v_j - v_i$ is orthogonal to $v_k - v_i$ for all $k \neq i$. Because $\Aff(S) = \R^{p, q}$, the vectors $v_k - v_i$ with $k \neq i$ span the whole vector space $\R^{p, q}$. Then, since the inner product on $\R^{p, q}$ is non-degenerate, we find that $v_i - v_j = 0$. However, this is a contradiction. Therefore $d_{i, j} > 0$ for any distinct $i, j$. Since $S$ is a non-obtuse set, all the triangle inequalities for $d_{i, j}$ are satisfied, and hence $d$ is indeed a metric on $X$. By the definition of $d_{i, j}$, we have
$$
E^{(0)} = 2 \mathrm{Gram}(v_1 - v_0, \ldots, v_{n-1} - v_0).
$$
Since $\mathrm{Aff}(S) = \R^{p, q}$, there exists a basis $\{ v_{i_1} - v_0, \ldots, v_{i_{p+q}} - v_0 \}$ of $\R^{p, q}$. From the basis, we can construct an orthonormal basis $\{ u_1, \ldots, u_{p+q} \}$ of $\R^{p, q}$. Then the linear expression of the vectors $v_i - v_0$ by the basis $\{ u_i \}$ provides us with an invertible matrix $A \in M(n-1, \R)$ such that
$$
\frac{1}{2}E^{(0)}
= {}^tA \ 
\mathrm{diag}(
\overbrace{1, \ldots, 1}^p, \overbrace{-1, \ldots, -1}^q, 0, \ldots, 0
) \
A.
$$
Hence the excess inertia of $d$ is $(p, q, n - p - q - 1)$.
\end{proof}

To illustrate Proposition \ref{prop:Schoenberg_embedding} and Proposition \ref{prop:metric_and_non_obtuse_set}, we give an example: In this paper, a $4$-point metric space will be called a \textit{metric rectangle} if its points can be cyclically ordered so that, for some $a, b > 0$, the four successive distances are $a, b, a, b$ and both diagonal distances are $a + b$. Then the distance matrix of a metric rectangle is described as
$$
D =
\left(
\begin{array}{cccc}
0 & a & a + b & b \\
a & 0 & b & a + b \\
a + b & b & 0 & a \\
b & a + b & a & 0
\end{array}
\right).
$$
One can realize this metric $d$ on $\{ 0, 1, 2, 3 \}$ as $d_{i, j} = \lVert v_i - v_j \rVert^2$, where $\{ v_0, v_1, v_2, v_3 \}$ is a non-obtuse set in the Euclidean space $\R^2$ given by
\begin{align*}
v_0 &= \left(\begin{array}{c} 0 \\ 0 \end{array}\right), &
v_1 &= \left(\begin{array}{c} \sqrt{a} \\ 0 \end{array}\right), &
v_2 &= \left(\begin{array}{c} \sqrt{a} \\ \sqrt{b} \end{array}\right), &
v_3 &= \left(\begin{array}{c} 0 \\ \sqrt{b} \end{array}\right).
\end{align*}
By this description, a metric rectangle has the excess inertia $(2, 0, 1)$.

\subsection{Non-obtuse sets and acute sets in Euclidean spaces}

For later use, we summarize some facts about non-obtuse and acute sets in Euclidean spaces.

\begin{prop} \label{prop:fact_non_obtuse_and_acute}
The following holds true.
\begin{enumerate}
\item[(a)]
(Danzer-Gr\"{u}nbaum \cite{DG})
A non-obtuse set in $\R^d$ can contain at most $2^d$ points. 

\item[(b)]
(Croft \cite{Croft})
An acute set in $\R^3$ can contain at most $5$ points.

\item[(c)]
(Croft \cite{Croft})
For any $6$-point non-obtuse set in $\R^3$, there exists no vertex at which all the angles are acute.

\item[(d)]
(Gerencs\'{e}r and Harangi \cite{GH})
There exists an acute set in $\R^d$ consisting of $2^{d-1} + 1$ points.

\end{enumerate}
\end{prop}

Some comments are in order: For (a), a $2^d$-point non-obtuse set in $\R^d$ is realized as $\{ (x_1, \ldots, x_d) \in \R^d |\ x_i \in \{ 0, 1 \} \}$ for example. For (c), let $S \subset \R^3$ be a $6$-point non-obtuse set. The meaning of (c) is that there is no $q \in S$ such that the angle $\angle p q r$ is acute for every pair of distinct points $p, r \in S \backslash \{ q \}$. Thus, for any $q \in S$, there are some other $p, r \in S$ such that $\angle p q r$ is a right angle. This result follows from Croft's classification of $6$-point non-obtuse sets in $\R^3$. The classification shows that any $6$-point non-obtuse set $S$ in $\R^3$ has one of the following forms (\cite{Croft}, pp.\ 305--306):
\begin{enumerate}
\item[($\alpha$)]
$S$ consists of the six vertices of a triangular prism.

\item[($\beta$)]
Consider a rectangle on a plane. Choose a pair of opposite vertices, and raise a copy of the pair of vertices by the same height. Then the two raised points and the four vertices of the original rectangle form $S$.

\item[($\gamma$)]
Consider a triangle on a plane, raise its copy, and rotate the raised copy through $\pi$ about its circumcenter. The three vertices of the raised triangle and those of the original triangle form $S$.

\end{enumerate}

\subsection{Some constructions}

Some constructions of metrics are summarized here. The first construction is known as the $1$-sum in \cite{DL}.

\begin{lem} \label{lem:1_sum}
From a metric space of excess inertia $(p, q, r)$, we can construct a metric
space with one additional point whose excess inertia is $(p+1, q, r)$.
\end{lem}

\begin{proof}
Let $(X, d)$ be a metric space with inertia $(p, q, r)$. Choose a base point $x_0 \in X$ and a real number $a > 0$. We adjoin a new point $y$ to $X$ to define $X_+ = X \cup \{ y \}$. We now extend the metric $d$ on $X$ to $X_+$ by setting
$$
d(y, x) = d(x, y) = a + d(x_0, x)
$$
for any $x \in X$. Then the  triangle inequalities on $X_+$ follow from those on $X$. If the $x_0$th excess matrix of $X$ is $E$, then the $x_0$th excess matrix of $X_+$ is $E \oplus (2a)$. 
\end{proof}

\begin{lem} 
\label{lem:diagonal}
If $E \in M(n-1, \R)$ is an excess matrix of an $n$-point metric space, then so is $E + \mathrm{diag}(a_1,\ldots,a_{n-1})$ for any $a_i \ge 0$.
\end{lem}

\begin{proof}
Let $d$ be a metric on $X = \{ 0, 1, \ldots, n - 1 \}$, and $E$ its $0$th excess matrix. Define a new metric $d'$ on $X$ by
$$
d'_{i, j} =
\left\{
\begin{array}{ll}
0, & (i = j) \\
d_{i, j} + \frac{a_i + a_j}{2}, & (i \neq j)
\end{array}
\right.
$$
where we put $a_0 = 0$. The triangle inequality for distinct $i, j, k$ is verified by
$$
d'_{i, k} + d'_{k, j} - d'_{i, j}
=
d_{i, k} + d_{k, j} - d_{i, j} + a_k \ge 0.
$$
The $(i, j)$-entry of the $0$th excess matrix of $d'$ is computed as
$$
d'_{0, i} + d'_{0, j} - d'_{i, j}
=
\left\{
\begin{array}{ll}
2 d_{0, i} + a_i, & (i = j) \\
d_{0, i} + d_{0, j} - d_{i, j}, & (i \neq j)
\end{array}
\right.
$$
and the proof is completed.
\end{proof}

\begin{lem}
\label{lem:two}
If there exists a metric on a finite set $X$ whose excess inertia is $(p, 2, 0)$, then there also exists a metric on $X$ whose excess inertia is $(p, 1, 1)$.
\end{lem}

\begin{proof}
Suppose that the two negative eigenvalues of an excess matrix $E$ of $X$ are different and $\lambda_1 < \lambda_2 < 0$. In this case, Lemma \ref{lem:diagonal} shows that $E - \lambda_2 I$ is an excess matrix of a metric on $X$, where $I$ is the identity matrix. The inertia of $E - \lambda_2 I$ is $(p, 1, 1)$.

Suppose then that the negative eigenvalue $\lambda < 0$ is of multiplicity two. In this case, choose an eigenvector $u = (u_i) \in \R^{p+2}$ corresponding to $\lambda$. Since $u \neq 0$, there is some index $i \in \{ 1, 2, \ldots, p+2 \}$ such that $u_i \neq 0$. Let $e_i \in \R^{p+2}$ be the vector whose $i$th entry is $1$ and all other entries are $0$. Because $E$ has no zero eigenvalue, the inertia of $E' = E + \epsilon e_i {}^te_i$ remains $(p, 2, 0)$ if $\epsilon > 0$ is sufficiently small. By Lemma \ref{lem:diagonal}, the matrix $E'$ is again an excess matrix of a metric on $X$. Now we are to prove that the two negative eigenvalues of $E'$ are distinct. Once this is established, the first part of the proof can be applied, and we get a metric on $X$ of excess inertia $(p, 1, 1)$. 

To prove that the two negative eigenvalues of $E'$ are distinct, it suffices to show that the kernel $\mathrm{Ker}(E' - \lambda I)$ of the linear map $E' - \lambda I$ is one-dimensional. Given a vector $v \in \R^{p+2}$ such that $(E' - \lambda I)v = 0$, we have
$$
0 = (E - \lambda I)v + \epsilon \langle e_i, v \rangle e_i.
$$
Taking the inner product with $u$ gives us
$$
0 = \langle (E - \lambda I)u, v \rangle + 
\epsilon \langle e_i, v \rangle \langle e_i, u \rangle
=
\epsilon \langle e_i, v \rangle \langle e_i, u \rangle.
$$
It follows that $\langle e_i, v \rangle = 0$. This allows us to express
$$
\mathrm{Ker}(E' - \lambda I)
=
\{ v \in \mathrm{Ker}(E - \lambda I) |\ \langle e_i, v \rangle = 0 \}.
$$
Since $\dim \mathrm{Ker}(E - \lambda I) = 2$ and $u \in \mathrm{Ker}(E - \lambda I)$ is such that $\langle e_i, u \rangle \neq 0$, we conclude that $\dim \mathrm{Ker}(E' - \lambda I) = 1$, and the proof is completed.
\end{proof}

\begin{lem}
\label{lem:acute}
Suppose that there is an $n$-point acute set in $\R^p$ whose affine span is $\R^p$. Then, for each $q = 0, 1, \ldots, n - p - 1$, there exists an $n$-point metric space whose excess inertia is $(p, q, n - p - q - 1)$. 
\end{lem}

\begin{proof}
The cases $n \le 2$ are immediate, so we assume $n \ge 3$. Let $\{ x_0, \ldots, x_{n-1} \} \subset \R^p$ be an acute set. Without loss of generality, we can assume $\Aff(\{ x_0, \ldots, x_p \}) = \R^p$. We let $m$ be
$$
m = \min \{
\langle x_i - x_j, x_k - x_j \rangle |\
\mbox{distinct $i, j, k \in \{ 0, 1, \ldots, n-1 \}$}
\}.
$$
This is positive, because $\{ x_0, \ldots, x_{n-1} \}$ is an acute set. For $i = p + 1, \ldots, p + q$, choose $t_i \in \R$ such that $0 < t_i^2 < m$. We then define $v_i \in \R^{p, q}$ for $i = 0, \ldots, n-1$ as follows. For $0 \le i \le p$ and $p + q < i \le n - 1$, we put $v_i = {}^t(x_i, 0) \in \R^{p, q}$. For $p + 1 \le i \le p + q$, we put
\begin{align*}
v_{p+1} &=
\left(
\begin{array}{c}
x_{p+1} \\ t_{p+1} \\ 0 \\ \vdots \\ 0
\end{array}
\right),
&
v_{p+2} &=
\left(
\begin{array}{c}
x_{p+2} \\ 0 \\ t_{p+2} \\ \vdots \\ 0
\end{array}
\right),
&
&\cdots
&
v_{p + q} &=
\left(
\begin{array}{c}
x_{p + q} \\ 0 \\ 0 \\ \vdots \\ t_{p + q}
\end{array}
\right).
\end{align*}
The affine span of $S = \{ v_0, \ldots, v_{n-1} \} \subset \R^{p, q}$ is clearly $\R^{p, q}$. Thus the proof will be completed by Proposition \ref{prop:metric_and_non_obtuse_set}, once $S$ is shown to be a non-obtuse set in $\R^{p, q}$. To show this, we let $I = \{ 0, 1, \ldots, p, p + q + 1, \ldots, n - 1 \}$ and $J = \{ p + 1, \ldots, p + q \}$. 
\begin{itemize}
\item
If $i, j, k \in I$ are distinct, then $\langle v_i - v_j, v_k - v_j \rangle_{p, q} = \langle x_i - x_j, x_k - x_j \rangle > 0$.

\item
If $i, j, k \in J$ are distinct, then $\langle v_i - v_j, v_k - v_j \rangle_{p, q} = \langle x_i - x_j, x_k - x_j \rangle - t_j^2 \ge m - t_j^2 > 0$.

\item
If distinct $i, j, k$ are such that $i, j \in I$ and $k \in J$, then
\begin{align*}
\langle v_i - v_j, v_i - v_k \rangle_{p, q}
&=
\langle x_i - x_j, x_i - x_k \rangle > 0,
\\
\langle v_j - v_i, v_j - v_k \rangle_{p, q}
&=
\langle x_j - x_i, x_j - x_k \rangle > 0,
\\
\langle v_k - v_i, v_k - v_j \rangle_{p, q}
&=
\langle x_k - x_i, x_k - x_j \rangle - t_k^2
\ge m - t_k^2 > 0.
\end{align*}

\item
If distinct $i, j, k$ are such that $i \in I$ and $j, k \in J$, then
\begin{align*}
\langle v_i - v_j, v_i - v_k \rangle_{p, q}
&=
\langle x_i - x_j, x_i - x_k \rangle > 0,
\\
\langle v_j - v_i, v_j - v_k \rangle_{p, q}
&=
\langle x_j - x_i, x_j - x_k \rangle - t_j^2 > 0,
\\
\langle v_k - v_i, v_k - v_j \rangle_{p, q}
&=
\langle x_k - x_i, x_k - x_j \rangle - t_k^2 > 0.
\end{align*}

\end{itemize}
Therefore $S$ is indeed an acute set in $\R^{p, q}$.
\end{proof}

\begin{lem} \label{lem:L1_product}
Let $(X, d)$ be an $n$-point metric space of excess inertia $(n_+, n_-, n_0)$, and $(X', d')$ an $n'$-point metric space of excess inertia $(n'_+, n'_-, n'_0)$. Then the $L^1$-product of $(X, d)$ and $(X', d')$ has the excess inertia 
$$
(n_+ + n'_+, n_- + n'_-, (n - 1)(n' -1) + n_0 + n'_0).
$$
\end{lem}

\begin{proof}
Recall that the $L^1$-product of $(X, d)$ and $(X', d')$ is a metric $\delta$ on the direct product $X \times X'$ given by
$$
\delta((x, x'), (y, y')) = d(x, y) + d'(x', y').
$$
In view of Proposition \ref{prop:Schoenberg_embedding}, suppose that $d$ on $X = \{ x_0, x_1, \ldots, x_{n-1} \}$ and $d'$ on $X' = \{ x'_0, \ldots, x'_{n'-1} \}$ are realized as 
\begin{align*}
d(x_i, x_j) &= \lVert u_i - u_j \rVert_{n_+, n_-}^2, &
d'(x'_i, x'_j) &= \lVert u'_i - u'_j \rVert_{n'_+, n'_-}^2, 
\end{align*}
in terms of sets of vectors $\{ u_0, u_1, \ldots, u_{n-1} \} \subset \R^{n_+, n_-}$ and $\{ u'_0, u'_1, \ldots, u'_{n' - 1} \} \subset \R^{n'_+, n'_-}$. Without loss of generality, we can assume that $u_0 = 0$ and $u'_0 = 0$. We can also assume that $\{ u_1, \ldots, u_{n_+ + n_-} \}$ and $\{ u'_1, \ldots, u'_{n'_+ + n'_-} \}$ generate $\R^{n_+, n_-}$ and $\R^{n'_+, n'_-}$, respectively. For $i = 0, 1, \ldots, n-1$ and $j = 0, 1, \ldots, n'-1$, we put 
$$
u_{i, j} = (u_i, u'_j) \in \R^{n_+, n_-} \times \R^{n'_+, n'_-}.
$$
The vector space $V = \R^{n_+, n_-} \times \R^{n'_+, n'_-}$ is equipped with the inner product
$$
\langle (v, v'), (w, w') \rangle_V
=
\langle v, w \rangle_{n_+, n_-} + \langle v', w' \rangle_{n'_+, n'_-}.
$$
It follows that 
$$
\lVert u_{i, j} - u_{k, \ell} \rVert^2_V
=
\lVert u_i - u_k \rVert^2_{n_+, n_-} +
\lVert u'_j - u'_\ell \rVert^2_{n'_+, n'_-}
=
\delta ((x_i, x'_j), (x_k, x'_\ell)).
$$
Hence the set of vectors $S = \{ u_{i, j} \}$ realizes the $L^1$-product $\delta$. The set $S$ spans $V$, because $u_{i, 0} = (u_i, 0)$ and $u_{0, j} = (0, u'_j)$ are contained in $S$. The inner product on $V$ has the inertia $(n_+ + n'_+, n_- + n'_-, 0)$. Therefore the numbers of positive and negative eigenvalues of the Gram matrix of the vectors in $S \backslash \{ u_{0, 0} \}$ are $n_+ + n'_+$ and $n_- + n'_-$, respectively. The number of zero eigenvalues is computed as
\begin{align*}
(nn' - 1) - (n_+ + n_-) - (n'_+ + n'_-)
&= (nn' - 1) - (n - 1 - n_0) - (n' - 1 - n'_0) \\
&= (n - 1)(n' - 1) + n_0 + n'_0.
\end{align*}
Hence the excess inertia of $\delta$ is as stated.
\end{proof}

\section{The classification of excess inertia}
\label{sec:classification_of_excess_inertia}

This section is the main part of this paper, and we study the excess inertia of finite metric spaces consisting of $n \le 10$ points. In the case of $n = 2, 3, 4$, the classifications are more or less known, but proofs are provided for the reader's convenience. The case of $n = 5$ was first established by Asao and Taira. Although our proof is slightly different, it is essentially the same as theirs. The cases of $6 \le n \le 10$ are our original results, where the case of $n = 6$ requires a non-trivial argument relying on Croft's classification of $6$-point non-obtuse sets in $\R^3$. In the case of $n = 10$, we only state a partial result.

\subsection{Two points and three points}

The cases of two and three points can be seen in an elementary way, and would be well known.

\begin{prop}  \label{prop:3_point}
The only excess inertia of a $2$-point metric space is $(1, 0, 0)$. The only excess inertia of a $3$-point metric space is $(2, 0, 0)$.
\end{prop}

\begin{proof}
In the $2$-point case, the excess matrix $E = (2d_{0, 1}) \in M(1, \R)$ is of inertia $(1, 0, 0)$ clearly. In the $3$-point case, the $0$th excess matrix is
$$
E =
\left(
\begin{array}{cc}
2 d_{0,1} & d_{0, 1} + d_{0, 2} - d_{1, 2} \\
d_{0, 1} + d_{0, 2} - d_{1, 2} & 2d_{0, 2}
\end{array}
\right).
$$
To simplify the notation, we put $a = d_{0, 1}$, $b = d_{0, 2}$, $c = d_{1, 2}$ and $u = (a + b - c)/2$. By the triangle inequalities, we have
\begin{align*}
a - u &= \frac{a + c - b}{2} \ge 0, &
b - u &= \frac{b + c - a}{2} \ge 0,
\end{align*}
and $u \ge 0$. It follows that $0 \le u^2 \le \min\{ a, b \}^2 \le ab$. Therefore $\det E = 4 (ab - u^2) \ge 0$. Suppose that $ab = u^2$. This implies $u = a = b$, and hence $c = a + b - 2u = 0$. This contradicts $c > 0$, and we conclude $\det E > 0$. Now, Sylvester's criterion implies that $E$ is positive definite and has the inertia $(2, 0, 0)$.
\end{proof}

\subsection{Four points}

The following classification of the excess inertia of a $4$-point metric space may also be well known to the experts, through the study of negative type metric spaces. 

\begin{prop} \label{prop:4_point}
The possible excess inertias of $4$-point metric spaces are precisely $(3, 0, 0)$ and $(2, 0, 1)$. In particular, the excess inertia of a $4$-point metric space is $(2, 0, 1)$ if and only if it is a metric rectangle.
\end{prop}

\begin{proof}
A proof based on techniques in metric geometry is given in \cite{HLMT} (Proposition 6.1, page 266) for example. We here give a direct proof.

We express the $0$th excess matrix of a metric $d$ on $X = \{ 0, 1, 2, 3 \}$ as
$$
E = 2
\left(
\begin{array}{ccc}
a & u & v \\
u & b & w \\
v & w & c
\end{array}
\right),
$$
where $a = d_{0, 1}$, $b = d_{0, 2}$, $c = d_{0, 3}$ and
\begin{align*}
u &= \frac{a + b - d_{1,2}}{2}, &
v &= \frac{a + c - d_{1,3}}{2}, &
w &= \frac{b + c - d_{2,3}}{2}.
\end{align*}
Renaming the points on $X$ if necessary, we can assume $u = \min\{ u, v, w \}$. This assumption implies $u \le v$ and $u \le w$. Then the triangle inequalities give $u \ge 0$ and
\begin{align*}
a - v &= \frac{d_{0, 1} + d_{1, 3} - d_{0, 3}}{2} \ge 0, \\
b - w &= \frac{d_{0, 2} + d_{2, 3} - d_{0, 3}}{2} \ge 0, \\
c + u - v - w &= \frac{d_{1, 3} + d_{2, 3} - d_{1, 2}}{2} \ge 0.
\end{align*}
For any $z = (z_i)_{i = 1}^3 \in \R^3$, one can directly verify the following expression of the quadratic form associated to $E/2$
\begin{align*}
\frac{1}{2}
\langle z, E z \rangle
&=
u (z_1 + z_2 + z_3)^2
+ (v - u)(z_1 + z_3)^2 + (w - u)(z_2 + z_3)^2
\\
&\quad
+ (a - v)z_1^2 + (b - w)z_2^2 + (c + u - v - w)z_3^2.
\end{align*}
All coefficients in this expression are non-negative, so $E$ is positive semidefinite. 

Now, suppose that $\det E = 0$. Then there exists $z = (z_i) \neq 0$ such that $Ez = 0$. It turns out that $z_i \neq 0$ for all $i$. To see this, let us assume $z_3 = 0$ for example. It follows that
$$
0 = \langle z, E z \rangle
=
( z_1 \ z_2) 
\left(
\begin{array}{cc}
2 d_{0, 1} & d_{0, 1} + d_{0, 2} - d_{1, 2} \\
d_{0, 1} + d_{0, 2} - d_{1, 2} & 2d_{0, 2}
\end{array}
\right)
\left(
\begin{array}{c}
z_1 \\ z_2 
\end{array}
\right).
$$
However, the excess matrix of the $3$-point subspace $\{ 0, 1, 2 \}$ is positive definite. So we get $z_1 = z_2 = 0$, which contradicts $z \neq 0$. In the same way, we find that $z_i \neq 0$ for all $i = 1, 2, 3$. This result and the expression of $\langle z, E z \rangle$ yield
\begin{align}
a &= v, &
b &= w, &
c + u &= v + w.
\tag{$\ast$} \label{formula:4_point_proof}
\end{align}
Assume here that $u > 0$. Then $z_1 + z_2 + z_3 = 0$. If $v > u$, then we have $z_1 + z_3 = 0$, and hence $z_2 = 0$. This is a contradiction, and we conclude $v = u$. In the same way, we find $w = u$. Together with (\ref{formula:4_point_proof}), we get $a = b = u$, but this results in a contradiction $d_{1, 2} = a + b - 2u = 0$. Therefore $u = 0$. In this case, (\ref{formula:4_point_proof}) gives rise to
\begin{align*}
d_{0, 1} &= d_{2, 3} = a, &
d_{0, 2} &= d_{1, 3} = b, &
d_{0, 3} &= d_{1, 2} = a + b,
\end{align*}
which is a metric rectangle, namely, a $4$-point metric space with excess inertia $(2, 0, 1)$. 

To conclude, the candidates for the excess inertia of a $4$-point metric space are $(3, 0, 0)$ and $(2, 0, 1)$. (The inertia $(1, 0, 2)$ does not appear.) The former inertia can be realized by applying the $1$-sum construction in Lemma \ref{lem:1_sum} to a $3$-point metric space, for example. The latter is realized by a metric rectangle. Finally, if the excess inertia of a $4$-point metric space $X$ is $(2, 0, 1)$, then $\det E = 0$, so that $X$ is a metric rectangle. 
\end{proof}

\subsection{Five points}

\begin{thm} \label{thm:5_point}
The possible excess inertias of $5$-point metric spaces are precisely the following.
$$
\begin{array}{|c|c|c|}
\hline
& n_0 = 0 & n_0 = 1 \\
\hline
n_- = 0 & (4, 0, 0) & (3, 0, 1) \\
\hline
n_- = 1 & (3, 1, 0) & \\
\hline
\end{array}
$$
\end{thm}

For the proof, we show a lemma.

\begin{lem} \label{lem:non_rectangle_in_5_point_metric_space}
Any $5$-point metric space contains at least two $4$-point subspaces which are not metric rectangles.
\end{lem}

\begin{proof}
In the case where all the $4$-point subspaces of $X = \{ 1, 2, 3, 4, 5 \}$ are not metric rectangles, we have nothing to prove. Thus, suppose that a subspace $Y_1 = \{ 2, 3, 4, 5 \} \subset X$ is a metric rectangle. We can assume that the distance matrix of $X$ is
$$
D_X
=
\left(
\begin{array}{ccccc}
0 & d_2 & d_3 & d_4 & d_5 \\
d_2 & 0 & a & a + b & b \\
d_3 & a & 0 & b & a + b \\
d_4 & a + b & b & 0 & a \\
d_5 & b & a + b & a & 0
\end{array}
\right),
$$
where $a, b > 0$. If $4$-point subspaces of $X$ other than $Y_1$ are not metric rectangles, then the claim of the lemma is verified. If there exists another $4$-point subspace which is a metric rectangle, then it shares three points with $Y_1$. Then, by renumbering the points if necessary, we can assume that $Y_5 = \{ 1, 2, 3, 4 \} \subset X$ is the other metric rectangle. The distance matrix of $Y_5$ is 
$$
D_{Y_5} =
\left(
\begin{array}{cccc}
0 & d_2 & d_3 & d_4 \\
d_2 & 0 & a & a + b \\
d_3 & a & 0 & b \\
d_4 & a + b & b & 0 \\
\end{array}
\right).
$$
The assumption that $Y_5$ is a metric rectangle yields
\begin{align*}
d_2 &= b, &
d_3 &= a+b, &
d_4 &= a.
\end{align*}
Here let us consider $Y_2 = \{ 1, 3, 4, 5 \}$, whose distance matrix is
$$
D_{Y_2}
=
\left(
\begin{array}{cccc}
0 & a + b & a & d_5 \\
a + b & 0 & b & a + b \\
a & b & 0 & a \\
d_5 & a + b & a & 0
\end{array}
\right).
$$
This subspace cannot be a metric rectangle for any choice of $a, b, d_5$. Indeed, for $Y_2$ to be a metric rectangle, $d_{3, 5} = a + b$ must be equal to $d_{1, 4} = a$. However, this is impossible. Similarly, the distance matrix of $Y_4 = \{ 1, 2, 3, 5 \}$ is

$$
D_{Y_4}
=
\left(
\begin{array}{cccc}
0 & b & a+b & d_5 \\
b & 0 & a & b \\
a+b & a & 0 & a + b \\
d_5 & b & a + b & 0
\end{array}
\right).
$$
For $Y_4$ to be a metric rectangle, $d_{3, 5} = a + b$ must be equal to $d_{1, 2} = b$, which is also impossible. Hence two $4$-point subspaces which are not metric rectangles are found, and the proof is completed. (Although this is unnecessary, the distance matrix of $Y_3 = \{ 1, 2, 4, 5 \}$ is
\begin{align*}
D_{Y_3}
&=
\left(
\begin{array}{cccc}
0 & b & a & d_5 \\
b & 0 & a + b & b \\
a & a + b & 0 & a \\
d_5 & b & a & 0
\end{array}
\right).
\end{align*}
This can be a rectangle if $d_5 = a + b$ and $a = b$.)
\end{proof}

\begin{proof}[Proof of Theorem \ref{thm:5_point}]
By Proposition \ref{prop:4_point}, the excess inertia of a $4$-point metric space is either $(3, 0, 0)$ or $(2, 0, 1)$. By Lemma \ref{lem:possible_inertia}, the excess inertia of a $5$-point metric space must be one of the following.$$
\begin{array}{|c|c|c|c|}
\hline
& n_0 = 0 & n_0 = 1 & n_0 = 2 \\
\hline
n_- = 0 & (4, 0, 0) & (3, 0, 1) & (2, 0, 2) \\
\hline
n_- = 1 & (3, 1, 0) & (2, 1, 1) & \\
\hline
\end{array}
$$
By Lemma \ref{lem:1_sum} and Proposition \ref{prop:4_point}, we can realize metric spaces of excess inertia $(4, 0, 0)$ and $(3, 0, 1)$. There exists no $5$-point metric space of excess inertia $(2, 0, 2)$. This is because the existence yields a $5$-point non-obtuse set in $\R^2$ by Proposition \ref{prop:metric_and_non_obtuse_set}, which contradicts the theorem of Danzer-Gr\"{u}nbaum as stated in Proposition \ref{prop:fact_non_obtuse_and_acute} (a). By Proposition \ref{prop:fact_non_obtuse_and_acute} (d), there exists a $5$-point acute set in $\R^3$. Thus, by Lemma \ref{lem:acute}, we get a $5$-point metric space of excess inertia $(3, 1, 0)$. Finally, there exists no $5$-point metric space of excess inertia $(2, 1, 1)$. If such a $5$-point metric space $X$ exists, then the Cauchy interlace theorem and Proposition \ref{prop:4_point} force the excess inertia of any $4$-point subspace in $X$ to be $(2, 0, 1)$. As is seen in Proposition \ref{prop:4_point}, a metric space of excess inertia $(2, 0, 1)$ is a metric rectangle. Thus, all the $4$-point subspaces in $X$ are metric rectangles. But, this contradicts Lemma \ref{lem:non_rectangle_in_5_point_metric_space}.
\end{proof}

\begin{rem} \label{lem:realize_(3,1,0)_by_K32}
An example of a $5$-point metric space of excess inertia $(3, 1, 0)$ is given by the complete bipartite graph $K_{3, 2}$, equipped with its shortest-path metric. The distance matrix and an excess matrix of $K_{3, 2}$ are
\begin{align*}
D &=
\left(
\begin{array}{ccccc}
0 & 2 & 1 & 1 & 1 \\
2 & 0 & 1 & 1 & 1 \\
1 & 1 & 0 & 2 & 2 \\
1 & 1 & 2 & 0 & 2 \\
1 & 1 & 2 & 2 & 0
\end{array}
\right),
&
E &=
2
\left(
\begin{array}{cccc}
2 & 1 & 1 & 1 \\
1 & 1 & 0 & 0 \\
1 & 0 & 1 & 0 \\
1 & 0 & 0 & 1
\end{array}
\right).
\end{align*}
It is easy to see that $E/2$ is the Gram matrix of $v_1, \ldots, v_4 \in \R^{3, 1}$ given by
\begin{align*}
v_1 &= \left(\begin{array}{c} 1 \\ 1 \\ 1 \\ 1 \end{array}\right), &
v_2 &= \left(\begin{array}{c} 1 \\ 0 \\ 0 \\ 0\end{array}\right), &
v_3 &= \left(\begin{array}{c} 0 \\ 1 \\ 0 \\ 0\end{array}\right), &
v_4 &= \left(\begin{array}{c} 0 \\ 0 \\ 1 \\ 0 \end{array}\right).
\end{align*}
Thus, putting $v_0 = 0$, we can realize the metric $d$ on $K_{3, 2}$ as $d_{i, j} = \lVert v_i - v_j \rVert^2_{3, 1}$. As a result, the excess inertia of $K_{3,2}$ is $(3, 1, 0)$.
\end{rem}

\subsection{Six points}

\begin{thm} \label{thm:6_point}
The possible excess inertias of $6$-point metric spaces are precisely the following.
$$
\begin{array}{|c|c|c|c|}
\hline
& n_0 = 0 & n_0 = 1 & n_0 = 2 \\
\hline
n_- = 0 & (5, 0, 0) & (4, 0, 1) & (3, 0, 2) \\
\hline
n_- = 1 & (4, 1, 0) & & \\
\hline
\end{array}
$$
\end{thm}

\begin{proof}
Using Lemma \ref{lem:possible_inertia} based on Theorem \ref{thm:5_point}, we can find the following candidates for the excess inertia of a $6$-point metric space.
$$
\begin{array}{|c|c|c|c|}
\hline
& n_0 = 0 & n_0 = 1 & n_0 = 2 \\
\hline
n_- = 0 & (5, 0, 0) & (4, 0, 1) & (3, 0, 2) \\
\hline
n_- = 1 & (4, 1, 0) & (3, 1, 1) & \\
\hline
n_- = 2 & (3, 2, 0) & & \\
\hline
\end{array}
$$
By Lemma \ref{lem:1_sum} and  Theorem \ref{thm:5_point}, we can realize metric spaces of excess inertia $(5, 0, 0)$, $(4, 0, 1)$ and $(4, 1, 0)$. In view of Proposition \ref{prop:metric_and_non_obtuse_set} and Proposition \ref{prop:fact_non_obtuse_and_acute} (a), we can construct a metric space of excess inertia $(3, 0, 2)$ from a $6$-point non-obtuse set in $\R^3$. We will prove the non-existence of a $6$-point metric space of excess inertia $(3, 1, 1)$ in Proposition \ref{prop:6_point_(3,1,1)}. Then, from Lemma \ref{lem:two}, it follows that there is no $6$-point metric space of excess inertia $(3, 2, 0)$.
\end{proof}

\begin{prop} \label{prop:6_point_(3,1,1)}
There exists no $6$-point metric space of excess inertia $(3, 1, 1)$.
\end{prop}

Since the proof of this proposition is relatively long and may appear to be complicated, we sketch the outline: First, we realize the hypothetical metric $d$ on $X = \{ 0, 1, \ldots, 5 \}$ of excess inertia $(3, 1, 1)$ as $d_{i, j} = \lVert v_i - v_j \rVert^2_{3, 1}$ by means of a non-obtuse set $S = \{ v_0, \ldots, v_5 \}$ in $\R^{3, 1}$. Then we deform the metric. For this purpose, we leave $S$ fixed, and consider a deformation of the inner product given by $\langle \ , \ \rangle_{3, 1} + \tau H ( \ , \ )$, where $\tau \ge 0$ is a deformation parameter, and $H$ is a positive definite inner product. The key point is that we construct $H$ such that $S \subset \R^4$ is non-obtuse with respect to $H$. Thanks to this property, $d$ is deformed within metrics. Then, at a particular value of the parameter, the deformed metric yields a $6$-point non-obtuse set in $\R^3$, which has a property contradicting Croft's result as stated in Proposition \ref{prop:fact_non_obtuse_and_acute} (c).

The inner product $H$ on $\R^4$ above will be constructed in Lemma \ref{lem:PD_inner_product} as the pull-back of the standard inner product on a larger space $\R^d$ under an injective linear map $T : \R^4 \to \R^d$. Note here that there exists an affine dependence relation of vectors $v_i$ in $S$, namely, a non-zero vector $\lambda = (\lambda_i)$ such that $\sum_i \lambda_i v_i = 0$ and $\sum_i \lambda_i = 0$. It is clear that $\sum_i \lambda_i T(v_i) = 0$. In view of this simple fact, $T$ is constructed first by choosing six points in $\R^d$ which have the same affine dependence relation as that of vectors in $S$, and then linearly extending the correspondence of six vectors in $S$ and six points in $\R^d$. In the present case, the affine dependence relation of vectors in $S$ can be seen as the sharing of a point by an $(r-1)$-simplex and an $(s-1)$-simplex whose vertices constitute $S$, where $(r, s) = (2, 4), (3, 3)$ (see Lemma \ref{lem:aff_indep_rel_continue} (a)). Then, based on a standard realization of an $(n-1)$-simplex in $\R^n$, we can realize an $(r-1)$-simplex and an $(s-1)$-simplex in $\R^r \times \R^s = \R^{r + s} = \R^6$ whose vertices have the same affine dependence relation with vectors in $S$ (see Lemma \ref{lem:supplemental_configuration_in_R6}). In our actual choice of six points in $\R^6$, the $(s-1)$-simplex is rescaled by $\sqrt{\rho} > 0$ to eventually make $S$ non-obtuse with respect to the pull-back $H$ of the inner product on $\R^d = \R^6$ under $T$.

We now start to prove lemmas for the proof of Proposition \ref{prop:6_point_(3,1,1)}.

\begin{lem} \label{lem:aff_indep_rel}
Suppose that there exists a non-obtuse set $S = \{ v_0, \ldots, v_5 \} \subset \R^{3, 1}$ such that $\Aff(S) = \R^{3, 1}$. 
\begin{itemize}
\item[(a)]
We define the space of affine dependence relations of vectors in $S$ by 
$$
K = \bigg\{ 
\lambda = (\lambda_i)_{i = 0}^5 \in \R^6 \bigg|\
\sum_{i = 0}^5 \lambda_i v_i = 0, \
\sum_{i = 0}^5 \lambda_i = 0
\bigg\}.
$$
Then this vector space $K$ is $1$-dimensional.

\item[(b)]
Any non-zero vector $\lambda = (\lambda_i)$ in $K$ has no zero entries: $\lambda_i \neq 0$. 

\end{itemize}
\end{lem}

\begin{proof}
(a) Since $\Aff(S) = \R^{3, 1}$, the vectors $\{ v_i - v_0 \}_{i = 1}^5$ span $\R^{3, 1}$. Hence there exists $(\lambda_i)_{i = 1}^5 \in \R^5$ such that $\sum_{i = 1}^5 \lambda_i (v_i - v_0) = 0$, and the vector space consisting of such vectors $(\lambda_i)_{i = 1}^5$ is one dimensional. Setting $\lambda_0 = - \sum_{i = 1}^5 \lambda_i$, we see $\mathrm{dim}K = 1$.

(b) Let $d$ be the metric on $X = \{ 0, 1, \ldots, 5 \}$ defined by $d_{i, j} = \lVert v_i - v_j \rVert^2_{3, 1}$. Suppose $\lambda_5 = 0$ for example. Then the dimension of the vector space $V \subset \R^{3, 1}$ spanned by $\{ v_i - v_0 \}_{i = 1}^4$ is at most three, and hence the rank of the excess matrix of $d$ on the $5$-point subspace $\{ 0, 1, 2, 3, 4 \}$ is at most three. Then the inertia of the excess matrix is $(3, 0, 1)$ according to the classification of the excess inertia of a $5$-point metric space in Theorem \ref{thm:5_point}. As a result, $\dim V = 3$, and the restriction of the inner product $\langle \ , \ \rangle_{3, 1}$ to $V$ is positive definite. Then, by transforming $S$ by an isometry of $\R^{3, 1}$, we can assume that $v_i = (x_i, 0)$ for $i = 0, \ldots, 4$ and $v_5 = (x_5, t_5)$, where $x_i \in \R^3$ and $t_5 \in \R$. This is a consequence of the so-called Witt extension theorem \cite{Ar}. We have $t_5 \neq 0$ by $\Aff(S) = \R^{3, 1}$. Now, consider the set $S' = \{ x_0, x_1, \ldots, x_5 \} \subset \R^3$. The points in $S'$ are distinct, since $\lVert x_i - x_j \rVert^2 = \lVert v_i - v_j \rVert^2_{3, 1} > 0$ for $0 \le i < j \le 4$ and $\lVert x_i - x_5 \rVert^2 = \lVert v_i - v_5 \rVert^2_{3, 1} + t_5^2 > 0$ for $0 \le i \le 4$. We can also see that $S'$ is non-obtuse: if $i, j, k \in \{ 0, 1, \ldots, 4 \}$ are distinct, then
$$
\langle x_i - x_j, x_k - x_j \rangle 
= \langle v_i - v_j, v_k - v_j \rangle_{3, 1} \ge 0.
$$
If $i, j \in \{ 0, 1, \ldots, 4 \}$ are distinct, then
\begin{align*}
\langle x_i - x_j, x_5 - x_j \rangle
&= \langle v_i - v_j, v_5 - v_j \rangle_{3, 1} \ge 0, 
\\
\langle x_i - x_5, x_j - x_5 \rangle
&= \langle v_i - v_5, v_j - v_5 \rangle_{3, 1} + t_5^2 > 0.
\end{align*}
In addition, all the angles $\angle x_i x_5 x_j$ are acute, as computed above. However, this contradicts Proposition \ref{prop:fact_non_obtuse_and_acute} (c) implied by Croft's classification of $6$-point non-obtuse sets in $\R^3$. Therefore we conclude that $\lambda_5 \neq 0$. The other entries $\lambda_i$ can also be handled in the same way.
\end{proof}

\begin{lem} \label{lem:aff_indep_rel_continue}
Suppose that there exists a $6$-point non-obtuse set $S \subset \R^{3, 1}$ such that $\Aff(S) = \R^{3, 1}$. 
\begin{itemize}
\item[(a)]
Renaming the vectors in $S$ and translating the set $S$ by an isometry of $\R^{3, 1}$, we can assume that $S$ consists of $\{ a_i \}_{i = 1}^r$ and $\{ b_j \}_{j = 1}^s$ with $(r, s) = (2, 4)$ or $(3, 3)$ such that there exist $\alpha_i, \beta_j \in \R$ satisfying
\begin{align*}
\sum_{i = 1}^r \alpha_i a_i &= \sum_{j = 1}^s \beta_j b_j = 0, &
\sum_{i = 1}^r \alpha_i &= \sum_{j = 1}^s \beta_j = 1, &
\alpha_i, \beta_j &> 0.
\end{align*}

\item[(b)]
We have $\alpha_i + \beta_j \le 1$ for all $i, j$.

\item[(c)]
If we put $\alpha_* = \max_i \alpha_i$ and $\beta_* = \max_j \beta_j$, then there exists $\rho > 0$ such that
$$
\frac{\beta_*}{1 - \beta_*} \le \rho \le \frac{1 - \alpha_*}{\alpha_*}.
$$
\end{itemize}
\end{lem}

\begin{proof}
By Lemma \ref{lem:aff_indep_rel} (a) and (b), there exist real numbers $\lambda_i \neq 0$ such that $\sum_{i = 0}^5 \lambda_i v_i = 0$ and $\sum_{i = 0}^5 \lambda_i = 0$. Let $r$ be the cardinality of $\{ i |\ \lambda_i > 0 \}$, and $s$ that of $\{ i |\ \lambda_i < 0 \}$. By design, $r + s = 6$. Changing the overall sign if necessary, we can assume $r \le s$. Now, renaming $v_i$ and $\lambda_i$, we can assume that $S$ consists of vectors $\{ a_i \}_{i = 1}^r$ and $\{ b_j \}_{j = 1}^s$ and the positive real numbers $\alpha_i, \beta_j > 0$ satisfy
\begin{align*}
\sum_{i = 1}^r \alpha_i a_i &= \sum_{j = 1}^s \beta_j b_j, &
\sum_{i = 1}^r \alpha_i &= \sum_{j = 1}^s \beta_j.
\end{align*}
By a translation on $\R^{3, 1}$, we can further assume that $\sum_{i = 1}^r \alpha_i a_i = \sum_{j = 1}^s \beta_j b_j = 0$. By a normalization, we can also assume that $\sum_{i = 1}^r \alpha_i = \sum_{j = 1}^s \beta_j = 1$. Hence (a) is established, up to the claim that $r \neq 1$ (and hence $(r, s) = (2, 4), (3, 3)$).

Under the above assumption about $a_i, b_j, \alpha_i, \beta_j$, we have
\begin{align*}
&
\sum_{\ell \neq j}\beta_\ell (b_\ell - a_i) -
\sum_{k \neq i}\alpha_k (a_k - b_j)
\\
&
=
\bigg( \sum_\ell \beta_\ell (b_\ell - a_i) \bigg) - \beta_j (b_j - a_i)
- \bigg( \sum_k \alpha_k (a_k - b_j) \bigg) + \alpha_i (a_i - b_j)
\\
&=
 - \sum_\ell \beta_\ell a_i + \sum_k \alpha_k b_j
+ (- \beta_j - \alpha_i)(b_j - a_i)
=
(1 - \beta_j - \alpha_i)(b_j - a_i)
\end{align*}
for each $i = 1, \ldots, r$ and $j = 1, \ldots, s$. Taking the inner product with $b_j - a_i$, we get
\begin{align*}
(1 - \alpha_i - \beta_j) \lVert b_j - a_i \rVert^2_{3, 1}
&=
\sum_{\ell \neq j} 
\beta_\ell \langle b_\ell - a_i, b_j - a_i \rangle_{3, 1}
+
\sum_{k \neq i}
\alpha_k \langle b_j - a_k, b_j - a_i \rangle_{3, 1}.
\end{align*}
This proves (b), namely, $\alpha_i + \beta_j \le 1$ for $i = 1, \ldots, r$ and $j = 1, \ldots, s$. Here, if we assume $r = 1$, then $\alpha_1 = 1$ and $\beta_j = 0$ for all $j = 1, \ldots, 5$. This contradicts $\sum_j \beta_j = 1$, and we conclude $r \neq 1$, the remaining claim in (a). 

Finally, since $\alpha_* > 0$, $\beta_* > 0$ and $\alpha_* + \beta_* \le 1$, we have $0 < \alpha_*, \beta_* \le 1$ and 
$$
\frac{1 - \alpha_*}{\alpha_*} - \frac{\beta_*}{1 - \beta_*}
=
\frac{1 - \alpha_* - \beta_*}{\alpha_*(1 - \beta_*)} \ge 0.
$$
This ensures the existence of $\rho > 0$ in (c).
\end{proof}

\begin{lem} \label{lem:supplemental_configuration_in_R6}
Suppose that there exists a $6$-point non-obtuse set $S \subset \R^{3, 1}$ such that $\Aff(S) = \R^{3, 1}$. Suppose further that $S$ consists of $a_1, \ldots, a_r, b_1, \ldots, b_s$ with $(r, s) = (2, 4)$ or $(3, 3)$ and there are $\alpha_i, \beta_j > 0$ satisfying the relations in Lemma \ref{lem:aff_indep_rel_continue} (a). Let $\rho > 0$ be as in Lemma \ref{lem:aff_indep_rel_continue} (c).
\begin{itemize}
\item[(a)]
There exists a $6$-point set $\mathcal{S} = \{ A_1, \ldots, A_r, B_1, \ldots, B_s \}$ in $\R^6 = \R^{6, 0}$ such that
\begin{align*}
\langle A_i, A_k \rangle 
&= \frac{\delta_{i, k}}{\alpha_i} - 1, 
&
\langle B_j, B_\ell \rangle 
&= 
\rho \bigg( \frac{\delta_{j, \ell}}{\beta_j} - 1 \bigg), &
\langle A_i, B_j \rangle &= 0
\end{align*}
for all $i, k \in \{ 1, \ldots, r \}$ and $j, \ell \in \{ 1, \ldots, s \}$. Furthermore, we have
$$
\sum_{i = 1}^r \alpha_i A_i 
=
\sum_{j = 1}^s \beta_j B_j = 0.
$$

\item[(b)]
The affine span of $\mathcal{S}$ in $\R^6$ is four dimensional, and the space of affine dependence relations of vectors in $\mathcal{S}$ is spanned by $(\alpha_1, \ldots, \alpha_r, -\beta_1, \ldots, -\beta_s)$.

\item[(c)]
The $6$-point set $\mathcal{S} \subset \R^6$ is non-obtuse.

\item[(d)]
If we choose $\rho$ appropriately, then there exists a point $q \in \mathcal{S}$ such that all the angles $\angle p q r$ at $q$ are acute.

\end{itemize}
\end{lem}

\begin{proof}
(a) We define $A_i, B_j$ in $\R^6 = \R^{r + s} = \R^r \times \R^s$ by
\begin{align*}
A_i &=
\big(
\bigg(
\frac{\delta_{i, p} - \alpha_p}{\sqrt{\alpha_p}}
\bigg)_{p = 1, \ldots, r}, 
\overbrace{0, \ldots, 0}^s
\big),
&
B_j &=
\big(
\overbrace{0, \ldots, 0}^r,
\sqrt{\rho}
\bigg(
\frac{\delta_{j, q} - \beta_q}{\sqrt{\beta_q}}
\bigg)_{q = 1, \ldots, s}
\big).
\end{align*}
Then the formulae for inner products are verified directly.

(b)
Given real numbers $c_1, \ldots, c_r$, the $p$th entry of $\sum_{i = 1}^r c_i A_i$ is
$$
\sum_{i = 1}^r c_i \frac{\delta_{i, p} - \alpha_p}{\sqrt{\alpha_p}}
=
\frac{c_p - \alpha_p \sum_i c_i}{\sqrt{\alpha_p}}.
$$
Hence $\sum_{i = 1}^r c_i A_i = 0$ is equivalent to $c_p = \alpha_p \sum_i c_i$ for all $p = 1, \ldots, r$. Similarly, $\sum_j d_j B_j = 0$ is equivalent to $d_q = \beta_q \sum_j d_j$ for all $q = 1, \ldots, s$. Using these equivalences, we can verify that $A_1 - B_s, \ldots, A_r - B_s, B_1 - B_s, \ldots, B_{s-1} - B_s$ span a subspace of $\R^6$ whose dimension is $(r + s - 1) - 1 = 4$. Thus $\mathrm{Aff}(\mathcal{S})$ is four dimensional. In the same way, we can see that the space of affine dependence relations of vectors in $\mathcal{S}$ is spanned by $(\alpha_1, \ldots, \alpha_r, - \beta_1, \ldots, - \beta_s)$.

(c) 
For distinct $i, j, k$, the formulae in (a) give
\begin{align*}
\langle A_i - A_j, A_k - A_j \rangle &= \frac{1}{\alpha_j} > 0, &
\langle B_i - B_j, B_k - B_j \rangle &= \frac{\rho}{\beta_j} > 0.
\end{align*}
For $i$ and distinct $j, k$, we have
\begin{align*}
\langle B_i - A_j, A_k - A_j \rangle &= \frac{1}{\alpha_j} > 0, &
\langle A_i - B_j, B_k - B_j \rangle &= \frac{\rho}{\beta_j} > 0.
\end{align*}
For $i$ and distinct $j, k$, we have
\begin{align*}
\langle A_i - B_j, A_i - B_k \rangle &= \frac{1 - \alpha_i}{\alpha_i} - \rho, &
\langle B_i - A_j, B_i - A_k \rangle 
&= \rho \bigg( \frac{1 - \beta_i}{\beta_i} \bigg) - 1.
\end{align*}
These are non-negative, because $\rho$ satisfies
$$
\frac{\beta_j}{1 - \beta_j} \le \frac{\beta_*}{1 - \beta_*}
\le \rho \le \frac{1 - \alpha_*}{\alpha_*} \le \frac{1 - \alpha_i}{\alpha_i}.
$$
Therefore $\mathcal{S} \subset \R^6$ is a non-obtuse set.

(d)
If $\alpha_* + \beta_* < 1$, then we can choose $\rho$ such that
$$
\frac{\beta_*}{1 - \beta_*} < \rho < \frac{1 - \alpha_*}{\alpha_*}.
$$
In this case, $\mathcal{S}$ is an acute set. If $\alpha_* + \beta_* = 1$, then 
$$
\rho = \frac{\beta_*}{1 - \beta_*} 
= \frac{1 - \alpha_*}{\alpha_*}.
$$
Now, the following cases can occur.
\begin{itemize}
\item
We have $\alpha_i < \alpha_*$ for some $i$. In this case, the angles at $A_i$ are acute. 

\item
We have $\beta_j < \beta_*$ for some $j$. In this case, the angles at $B_j$ are acute. 

\item
We have $\alpha_i = \alpha_*$ for all $i$ and $\beta_j = \beta_*$ for all $j$. In this case, $\alpha_i = 1/r$ and $\beta_j = 1/s$. This yields a contradiction
$$
1 = \alpha_* + \beta_*
= \frac{1}{r} + \frac{1}{s}
=
\left\{
\begin{array}{ll}
\frac{3}{4}, ((r, s) = (2, 4)) \\
\frac{2}{3}, ((r, s) = (3, 3))
\end{array}
\right.
$$
so this case does not occur.

\end{itemize}
\end{proof}

\begin{lem} \label{lem:PD_inner_product}
Suppose that there exists a $6$-point non-obtuse set $S \subset \R^{3, 1}$ such that $\Aff(S) = \R^{3, 1}$. Then there exists a positive definite inner product $H = H( \ , \ )$ on $\R^4 = \R^{3, 1}$ satisfying the following properties.
\begin{itemize}
\item[(a)]
The set $S \subset \R^4$ is non-obtuse with respect to $H$.

\item[(b)]
There exists $q \in S$ such that all the angles $\angle p q r$ with respect to $H$ are acute.

\end{itemize}
\end{lem}

\begin{proof}
We can assume that $S$ consists of $a_1, \ldots, a_r, b_1, \ldots, b_s$ with $(r, s) = (2, 4)$ or $(3, 3)$ and there are $\alpha_i, \beta_j > 0$ satisfying the relations in Lemma \ref{lem:aff_indep_rel_continue} (a). Then, an appropriate choice of $\rho > 0$ gives us a non-obtuse set $\mathcal{S}$ in $\R^6$ consisting of $A_1, \ldots, A_r$, $B_1, \ldots, B_s$ such as in Lemma \ref{lem:supplemental_configuration_in_R6}. 

In this setup, we construct a linear map $T : \R^4 \to \R^6$. For its construction, recall that $\Aff(S) = \R^{3, 1} = \R^4$. Thus, for any $v \in \R^4$, there are $c_i$ and $d_j$ such that $v = \sum_{i = 1}^r c_i a_i + \sum_{j = 1}^s d_j b_j$ and $\sum_{i = 1}^r c_i + \sum_{j = 1}^s d_j = 1$. We then define $T : \R^4 \to \R^6$ by
$$
T(v) = \sum_{i = 1}^r c_i A_i + \sum_{j = 1}^s d_j B_j.
$$
The expression of $v$ above is not unique. If we have another expression $v = \sum_{i = 1}^r c'_i a_i + \sum_{j = 1}^s d'_j b_j$ in terms of $c'_i, d'_j$ such that $\sum_{i = 1}^r c'_i + \sum_{j = 1}^s d'_j = 1$, then there exists a unique $c \in \R$ such that $c_i - c'_i = c \alpha_i$ for all $i$ and $d_j - d'_j = - c \beta_j$ for all $j$, as a consequence of Lemma \ref{lem:aff_indep_rel} (a). Then Lemma \ref{lem:supplemental_configuration_in_R6} (a) allows us to conclude that $T$ is well-defined. By Lemma \ref{lem:supplemental_configuration_in_R6} (b), the affine span of $\mathcal{S}$ is four dimensional. It follows that $T : \R^4 \to \R^6$ is injective. 

Now, for $u, v \in \R^4$, we put $H(u, v) = \langle T(u), T(v) \rangle$ to define a symmetric bilinear form. Since $T$ is injective, $H$ is positive definite. By definition, $T(a_i) = A_i$ and $T(b_j) = B_j$. We then get
$$
H(a_i - a_j, a_k - a_j)
=
\langle T(a_i - a_j), T(a_k - a_j) \rangle
=
\langle A_i - A_j, A_k - A_j \rangle
$$
for distinct $i, j, k$. Thus, Lemma \ref{lem:supplemental_configuration_in_R6} (c) and (d) respectively imply (a) and (b) in this lemma.
\end{proof}

\begin{lem} \label{lem:contradicting_non_obtuse_set_in_R3}
Suppose that there exists a $6$-point non-obtuse set $S \subset \R^{3, 1}$ such that $\Aff(S) = \R^{3, 1}$. Then there exists a $6$-point non-obtuse set $\overline{S} \subset \R^3$ which admits a point $\bar{q} \in \overline{S}$ such that all the angles $\angle \bar{p} \bar{q} \bar{r}$ at $\bar{q}$ are acute.
\end{lem}

\begin{proof}
In this proof, we tentatively put $L(u, v) = \langle u, v \rangle_{3, 1}$ for $u, v \in \R^4 = \R^{3, 1}$. Let $H$ be an inner product in Lemma \ref{lem:PD_inner_product}. Since $H$ is positive definite, there is an orthonormal basis $\{ u_i \}_{i = 1}^4$ of $\R^4$ with respect to $H$. Then, diagonalizing the symmetric matrix $(L(u_i, u_j))_{i, j = 1}^4$, we get another orthonormal basis $\{ w_i \}_{i = 1}^4$ with respect to $H$ for which $L$ is diagonalized. Consequently, if we define a linear transformation $\phi : \R^4 \to \R^4$ by $\phi(x) = \sum_i x_i w_i$ for $x = (x_i) \in \R^4$, then the inner products (bilinear forms) $\phi^*H$ and $\phi^*L$ are expressed as
\begin{align*}
(\phi^*H)(x, y) &= \sum_{i = 1}^4 x_iy_i, &
(\phi^*L)(x, y) &= \sum_{i = 1}^4 \lambda_ix_iy_i,
\end{align*}
where $\lambda_i$ are the eigenvalues of $(L(u_i, u_j))_{i, j = 1}^4$. Because $(L(u_i, u_j))_{i, j = 1}^4$ has the inertia $(3, 1, 0)$, we can assume $\lambda_1 < 0 < \lambda_2 \le \lambda_3 \le \lambda_4$ without loss of generality. 

For $\tau > 0$, we define a symmetric bilinear form $L_\tau$ on $\R^4$ by
$$
L_\tau(u, v) = L(u, v) + \tau H(u, v).
$$
The pull-back inner product $\phi^*L_\tau$ has the expression
$$
(\phi^*L_\tau)(x, y) = \sum_{i = 1}^4 (\lambda_i + \tau) x_iy_i.
$$
We then choose $\tau$ to be $\tau_* = - \lambda_1 > 0$. Then $L_{\tau^*}$ has the inertia $(3, 0, 1)$. It is clear that an eigenvector of $\phi^*L_{\tau^*}$ whose eigenvalue is zero is $e_1 = (1, 0, 0, 0)$. If we define a linear surjection $\Pi : \R^4 \to \R^3$ by
$$
\Pi(x) = (
\sqrt{\lambda_2 + \tau_*}x_2, 
\sqrt{\lambda_3 + \tau_*}x_3,
\sqrt{\lambda_4 + \tau_*}x_4
)
$$
for $x = (x_i) \in \R^4$, then we apparently have
$$
\langle \Pi(x), \Pi(y) \rangle = (\phi^*L_{\tau^*})(x, y)
$$
for $x, y \in \R^4$. 

Now, we define $\overline{S} \subset \R^3$ by the image $\overline{S} = \Pi(\phi^{-1}(S))$ of $S \subset \R^4$ under the linear surjection $\Pi \circ \phi^{-1} : \R^4 \to \R^3$. For $s \in S$, we put $\bar{s} = \Pi(\phi^{-1}(s))$. For distinct $s, s' \in S$, we have
\begin{align*}
\lVert \bar{s} - \bar{s}' \rVert^2
=
\lVert \Pi(\phi^{-1}(s - s')) \rVert^2
&=
L_{\tau^*}(s - s', s - s')
\\
&=
L(s - s', s - s') + \tau_* H(s - s', s - s') > 0.
\end{align*}
Thus, $\overline{S}$ consists of six points. For distinct $p, q, r \in S$, we also have
$$
\langle \bar{p} - \bar{q}, \bar{r} - \bar{q} \rangle
=
L(p - q, r - q)
+ \tau_* H(p - q, r - q) \ge 0,
$$
since $S \subset \R^4$ is non-obtuse with respect to both $L$ and $H$. By Lemma \ref{lem:PD_inner_product}, there exists $q \in S$ at which all the angles $\angle p q r$ with respect to $H$ are acute. Then the computation above shows that all the angles $\angle \bar{p} \bar{q} \bar{r}$ at $\bar{q}$ are acute too. 
\end{proof}

\begin{proof}[Proof of Proposition \ref{prop:6_point_(3,1,1)}]
If there is a $6$-point metric space whose excess inertia is $(3, 1, 1)$, then we have a $6$-point non-obtuse set $S$ in $\R^{3, 1}$ such that $\Aff(S) = \R^{3, 1}$. Then, by Lemma \ref{lem:contradicting_non_obtuse_set_in_R3}, there exists a $6$-point non-obtuse set $\overline{S}$ in $\R^3$ equipped with a point $\bar{q}$ at which all the angles $\angle \bar{p}\bar{q}\bar{r}$ are acute. However, this contradicts a consequence of Croft's classification of $6$-point non-obtuse sets in $\R^3$ stated as Proposition \ref{prop:fact_non_obtuse_and_acute} (c).
\end{proof}

\smallskip

\begin{rem}
From Croft's classification \cite{Croft}, it also follows that any $6$-point non-obtuse set in $\R^3$ contains the four vertices of a Euclidean rectangle. It is possible to formulate a slightly different proof of Proposition \ref{prop:6_point_(3,1,1)} based on this fact. 
\end{rem}

\subsection{Seven points}

\begin{thm} \label{thm:7_point}
The possible excess inertias of $7$-point metric spaces are precisely the following.
$$
\begin{array}{|c|c|c|c|c|}
\hline
& n_0 = 0 & n_0 = 1 & n_0 = 2 & n_0 = 3 \\
\hline
n_- = 0 & (6, 0, 0) & (5, 0, 1) & (4, 0, 2) & (3, 0, 3)\\
\hline
n_- = 1 & (5, 1, 0) & (4, 1, 1) & &  \\
\hline
n_- = 2 & (4, 2, 0) & & & \\
\hline
\end{array}
$$
\end{thm}

\begin{proof}
By Lemma \ref{lem:possible_inertia}, the candidates for the excess inertia of a $7$-point metric space are those listed in the present theorem and $(3, 1, 2)$. However, this inertia is impossible: If such a $7$-point metric space $X$ exists, then the Cauchy interlace theorem and the classification of the $6$-point excess inertia in Theorem \ref{thm:6_point} imply that all the principal $5 \times 5$ submatrices of an excess matrix $E$ of $X$ have the same inertia $(3, 0, 2)$. We have $\det E = 0$. Thus, Sylvester's criterion implies that $E$ is positive semidefinite, which is a contradiction.

The $1$-sum construction in Lemma \ref{lem:1_sum} yields metric spaces of excess inertia $(6, 0, 0)$, $(5, 0, 1)$, $(4, 0, 2)$ and $(5, 1, 0)$. A $7$-point non-obtuse set in $\R^3$ yields a metric space of excess inertia $(3, 0, 3)$. By Proposition \ref{prop:fact_non_obtuse_and_acute} (d), there exists a $7$-point acute set in $\R^4$. Its affine span is $\R^4$: Otherwise, we get a $7$-point acute set in $\R^d$ with $d \le 3$, which contradicts Proposition \ref{prop:fact_non_obtuse_and_acute} (b). From this $7$-point acute set in $\R^4$, we get metric spaces of excess inertia $(4, 2, 0)$ and $(4, 1, 1)$.
\end{proof}

\begin{rem}
One can realize the inertia in Theorem \ref{thm:7_point} by a graph. For example, the complement $\overline{C_7}$ of the cycle graph $C_7$ realizes $(4, 2, 0)$, and the complete tripartite graph $K_{3, 2, 2}$ realizes $(4, 1, 1)$.
\end{rem}

\subsection{Eight points}

\begin{thm} \label{thm:8_point}
The possible excess inertias of $8$-point metric spaces are precisely the following.
$$
\begin{array}{|c|c|c|c|c|c|}
\hline
& n_0 = 0 & n_0 = 1 & n_0 = 2 & n_0 = 3 & n_0 = 4 \\
\hline
n_- = 0 & (7, 0, 0) & (6, 0, 1) & (5, 0, 2) & (4, 0, 3) & (3, 0, 4) \\
\hline
n_- = 1 & (6, 1, 0) & (5, 1, 1) & (4, 1, 2) & & \\
\hline
n_- = 2 & (5, 2, 0) & (4, 2, 1) & & & \\
\hline
n_- = 3 & (4, 3, 0) & & & & \\
\hline
\end{array}
$$
\end{thm}

\begin{proof}
By Lemma \ref{lem:possible_inertia}, the candidates for the excess inertia of an $8$-point metric space are those listed in the present theorem and $(3, 1, 3)$. However, $(3, 1, 3)$ is impossible. If an $8$-point metric space $X$ has the excess inertia $(3, 1, 3)$, then the Cauchy interlace theorem and Theorem \ref{thm:7_point} imply that all the principal $6 \times 6$ submatrices of an excess matrix $E$ of $X$ have the same inertia $(3, 0, 3)$. Because $\det E = 0$, Sylvester's criterion implies that $E$ is positive semidefinite, which is a contradiction.

The possible $11$ types of inertia are realized as follows:
\begin{itemize}
\item
The $1$-sum construction in Lemma \ref{lem:1_sum} gives $(7, 0, 0)$, $(6, 0, 1)$, $(5, 0, 2)$, $(4, 0, 3)$, $(6, 1, 0)$, $(5, 1, 1)$, $(5, 2, 0)$.

\item
From an $8$-point non-obtuse set in $\R^3$, we get $(3, 0, 4)$.

\item
By Proposition \ref{prop:fact_non_obtuse_and_acute}, there is an $8$-point acute set in $\R^4$, and its affine span is $\R^4$. Then, by Lemma \ref{lem:acute}, we get $(4, 3, 0)$, $(4, 2, 1)$, $(4, 1, 2)$.

\end{itemize}
This completes the proof.
\end{proof}

\subsection{Nine points}

\begin{thm} \label{thm:9_point}
The possible excess inertias of $9$-point metric spaces are precisely the following.
$$
\begin{array}{|c|c|c|c|c|c|c|c|}
\hline
& n_0 = 0 & n_0 = 1 & n_0 = 2 & n_0 = 3 & n_0 = 4  \\
\hline
n_- = 0 & (8, 0, 0) & (7, 0, 1) & (6, 0, 2) & (5, 0, 3) & (4, 0, 4) \\
\hline
n_- = 1 & (7, 1, 0) & (6, 1, 1) & (5, 1, 2) & (4, 1, 3) &  \\
\hline
n_- = 2 & (6, 2, 0) & (5, 2, 1) & (4, 2, 2) & & \\
\hline
n_- = 3 & (5, 3, 0) & (4, 3, 1) &  & & \\
\hline
n_- = 4 & (4, 4, 0)&  & & & \\
\hline
\end{array}
$$
\end{thm}

\begin{proof}
By Lemma \ref{lem:possible_inertia}, the candidates for the excess inertia of a $9$-point metric space are those listed in the present theorem together with $(3, 0, 5)$ and $(3, 1, 4)$. The inertia $(3, 0, 5)$ is impossible, because its existence implies that of a $9$-point non-obtuse set in $\R^3$. The inertia $(3, 1, 4)$ is also impossible. If we assume the existence of a metric space $X$ of excess inertia $(3, 1, 4)$, then the Cauchy interlace theorem, Theorem \ref{thm:8_point} and Sylvester's criterion imply that an excess matrix $E$ of $X$ is positive semidefinite, a contradiction. 

The possible $15$ types of inertia are realized as follows:
\begin{itemize}
\item
The $1$-sum constructs $(8, 0, 0)$, $(7, 0, 1)$, $(6, 0, 2)$, $(5, 0, 3)$, $(4, 0, 4)$, $(7, 1, 0)$, $(6, 1, 1)$, $(5, 1, 2)$, $(6, 2, 0)$, $(5, 2, 1)$, $(5, 3, 0)$.

\item
By Proposition \ref{prop:fact_non_obtuse_and_acute}, there is a $9$-point acute set in $\R^4$, and its affine span is $\R^4$. Then Lemma \ref{lem:acute} yields $(4, 4, 0)$, $(4, 3, 1)$, $(4, 2, 2)$, $(4, 1, 3)$.

\end{itemize}
This completes the proof.
\end{proof}

\subsection{Ten points}

In the $10$-point case, the classification of the excess inertia is not yet completed. A partial answer is the following.

\begin{prop} \label{prop:10_point}
Each inertia in the following list can be realized as the excess inertia of a $10$-point metric space.
$$
\begin{array}{|c|c|c|c|c|c|c|c|}
\hline
& n_0 = 0 & n_0 = 1 & n_0 = 2 & n_0 = 3 & n_0 = 4 & n_0 = 5 \\
\hline
n_- = 0 & (9,0,0) & (8,0,1) & (7,0,2) & (6,0,3) & (5,0,4) & (4,0,5) \\
\hline
n_- = 1 & (8,1,0) & (7,1,1) & (6,1,2) & (5,1,3) & (4,1,4) & \\
\hline
n_- = 2 & (7,2,0) & (6,2,1) & (5,2,2) & & & \\
\hline
n_- = 3 & (6,3,0) & (5,3,1) & & & & \\
\hline
n_- = 4 & (5,4,0) & & & & & \\
\hline
\end{array}
$$
\end{prop}

\begin{proof}
By Lemma \ref{lem:possible_inertia}, the candidates for the excess inertia of a $10$-point metric space are those listed in the present proposition together with $(4, q, 5 - q)$ for $q = 2, 3, 4, 5$. The inertias with $n_+ \ge 5$ can all be realized by the $1$-sum construction. The inertia $(4, 0, 5)$ can be realized by a $10$-point non-obtuse set in $\R^4$. Finally, the inertia $(4, 1, 4)$ is realized by the $L^1$-product of $K_{3, 2}$ and $K_2$. Here $K_{3, 2}$ stands for the bipartite graph $K_{3, 2}$ regarded as a $5$-point metric space by the shortest path length. As noted in Remark \ref{lem:realize_(3,1,0)_by_K32}, its excess inertia is $(3, 1, 0)$. Then $K_2$ stands for the $2$-point metric space such that the distance of the two points is $1$. By Lemma \ref{lem:L1_product}, the excess inertia of the $L^1$-product is computed as
$$
(3 + 1, 1 + 0, (5-1)(2-1) + 0 + 0) = (4, 1, 4),
$$
and the proof is completed.
\end{proof}

The above argument leaves open the realizability of the inertia $(4, q, 5-q)$ for $q = 2, 3, 4, 5$. If there exists a $10$-point acute set in $\R^4$, then all these inertias can be realized at once by Lemma \ref{lem:acute}. However, this problem also seems to be open.

\section{Applications}
\label{sec:application}

This section presents the two applications mentioned in \S\ref{sec:introduction}.

\subsection{The maximal number of points in non-obtuse sets}

\begin{lem} \label{lem:finiteness}
For any non-negative integers $p, q$, suppose that a set $S \subset \R^{p, q}$ satisfies the following three properties.
\begin{enumerate}
\item[(i)]
We have $\lVert u - v \rVert^2_{p, q} > 0$ for distinct $u, v \in S$.

\item[(ii)]
We have $\langle u - v, w - v \rangle_{p, q} \ge 0$ for distinct $u, v, w \in S$.

\item[(iii)]
We have $\Aff(S) = \R^{p, q}$

\end{enumerate}
If $S \neq \emptyset$, then $S$ consists of a finite number of points and its cardinality $\# S$ satisfies $p + q + 1 \le \# S \le 2^{p+q}$.
 \end{lem}

\begin{proof}
In the case where $S$ consists of one point, the third property of $S$ implies $p = q = 0$. In this case, the first and second properties are vacuously satisfied. We have $p + q + 1 \le \# S \le 2^{p+q}$ trivially. In the case where $S$ consists of two points, the third property of $S$ implies $(p, q) = (1, 0)$ or $(0, 1)$. The first property implies $(p, q) = (1, 0)$. The second property is vacuously satisfied. Then $p + q + 1 \le \# S \le 2^{p+q}$ trivially again. 

Henceforth we assume that $S$ contains at least three points. Then we apply the proof of Danzer and Gr\"{u}nbaum \cite{AZ,DG} to our setup: In view of the third property $\Aff(S) = \R^{p, q}$ of $S$, we can choose a finite subset $S' \subset S$ such that $\Aff(S') = \R^{p, q}$. If we put $n = \# S'$, then $n \ge p + q + 1$. Let $P$ be the convex hull of $S' = \{ v_1, \ldots, v_n \}$
$$
P =
\bigg\{
\sum_{i = 1}^n \lambda_i v_i \in \R^{p, q} 
\bigg|\
\sum_{i = 1}^n \lambda_i = 1, 
\lambda_i \ge 0 \ (i = 1, \ldots, n)
\bigg\}.
$$
For $v \in S'$, we let $f_v : \R^{p, q} \to \R^{p, q}$ be $f_v(x) = v + (x - v)/2 = (v + x)/2$. Using $f_v$, we define $P_v = f_v(P)$ for $v \in S'$. Then $P_v$ has the following properties.
\begin{itemize}
\item[(i)]
We have $P_v \subset P$ for any $v \in S'$. To see this, we express any $x \in P_v$ as $x = f_v(y)$ by using $y \in P$. We can express $y \in P$ as $y = \sum_i \lambda_i v_i$ by using $\lambda_i \ge 0$ such that $\sum_i \lambda_i = 1$. Let $j \in \{ 1, \ldots, n \}$ be the index such that $v_j = v$. Then we have
$$
x = 
\frac{v + y}{2}
=
\frac{1 + \lambda_j}{2} v_j
+ \sum_{i \neq j} \frac{\lambda_i}{2}v_i
\in P,
$$
because $\frac{1 + \lambda_j}{2} + \sum_{i \neq j} \frac{\lambda_i}{2} = 1$.

\item[(ii)]
For any distinct $u, v \in S'$, the intersection of the interiors of $P_u$ and $P_v$ in $\R^{p, q} = \R^{p + q}$ is empty
$$
\mathrm{Int}(P_u) \cap \mathrm{Int}(P_v) = \emptyset.
$$
To see this, suppose that there exists a point $z \in \mathrm{Int}(P_u) \cap \mathrm{Int}(P_v)$. Then there exist $x, y \in P$ such that
\begin{align*}
z &= f_u(x) = \frac{u + x}{2}, &
z &= f_v(y) = \frac{v + y}{2}.
\end{align*}
Thus, we get $x - y = v - u$, and 
\begin{align*}
\langle x, v - u \rangle_{p, q} - \langle y, v - u \rangle_{p, q}
=
\langle x - y, v - u \rangle_{p, q}
=
\lVert v  - u \rVert^2_{p, q}.
\end{align*}
We can then compute this in another way to see a contradiction: Since $z \in \mathrm{Int}(P_u)$, we also have $x \in \mathrm{Int}(P)$. Hence $x$ is expressed as $x = \sum_i \lambda_i v_i$ in terms of $\lambda_i > 0$ such that $\sum_i \lambda_i = 1$. This expression gives us
\begin{align*}
\langle v, v - u \rangle_{p, q} - \langle x, v - u \rangle_{p, q}
&=
\langle v - x, v - u \rangle_{p, q}
\\
&=
\langle \sum_i \lambda_i v - \sum_i \lambda_i v_i, 
v - u \rangle_{p, q}
\\
&=
\sum_i \lambda_i \langle v - v_i, v - u \rangle_{p, q}.
\end{align*}
This is non-negative, since $S' \subset S \subset \R^{p, q}$ is non-obtuse (the second property of $S$). We can further show that this is positive: If $\sum_i \lambda_i \langle v - v_i, v - u \rangle_{p, q} = 0$, then $\langle v - v_i, v - u \rangle_{p, q} = 0$ for all $i$. This means that $v - u$ is orthogonal to $v - v_i$ for all $i$. Note that $v = v_j$ for some $j \in \{ 1, \ldots, n \}$ and the span of $\{ v_j - v_i \}_{i \neq j}$ is the whole space $\R^{p, q}$ by the assumption $\Aff(S') = \R^{p, q}$. Because the inner product on $\R^{p, q}$ is non-degenerate, it follows that $v - u = 0$, which is a contradiction. As a consequence, we see
$$
\langle v, v - u \rangle_{p, q} - \langle x, v - u \rangle_{p, q} > 0.
$$
Similarly, we see $\langle y, v - u \rangle_{p, q} - \langle u, v - u \rangle_{p, q} > 0$. Hence we have
$$
\lVert v - u \rVert^2_{p, q}
=
\langle v, v - u \rangle_{p, q} - \langle u, v - u \rangle_{p, q} 
>
\langle x, v - u \rangle_{p, q} - \langle y, v - u \rangle_{p, q}.
$$
However, this contradicts the previous computation.

\end{itemize}

Now, let $\mu$ be the standard Lebesgue measure on $\R^{p, q} = \R^{p+q}$. By $\Aff(S') = \R^{p+q}$, we have $\mu(P) > 0$. Recalling the definition of $f_v$, we see $\mu(P_v) = \mu(P)/2^{p+q}$. Then, by the properties (i) and (ii) of $P_v$, we get
$$
\mu(P) \ge \mu(\cup_{i = 1}^n P_{v_i}) = \sum_{i = 1}^n \mu(P_{v_i})
= \frac{n}{2^{p+q}} \mu(P),
$$
and hence $\# S' = n \le 2^{p+q}$. As a result, $S$ is a finite set and $\# S \le 2^{p+q}$: Otherwise, we can enlarge $S'$ to a finite subset of $S$ consisting of $1 + 2^{p+q}$ points, which is a contradiction. Finally, by the third property $\Aff(S) = \R^{p, q} = \R^{p+q}$ of $S$, it must contain at least $p + q + 1$ points. 
\end{proof}

\begin{dfn}
For any non-negative integers $p, q$, we define $g(p, q)$ to be the maximal number of points contained in such a set $S \subset \R^{p, q}$ as in Lemma \ref{lem:finiteness}. We set $g(p, q) = 0$ if no such set exists.
\end{dfn}

\begin{lem} \label{lem:gpq_general_property}
The following holds true.
\begin{itemize}
\item[(a)]
For any non-negative integers $p, q$, we have either $g(p, q) = 0$ or
$$
p + q + 1 \le g(p, q) \le 2^{p+q}.
$$

\item[(b)]
For any non-negative integers $p, q$, the number $g(p, q)$ agrees with the maximal number $n$ such that there exists an $n$-point metric space whose excess inertia is $(p, q, n - p - q - 1)$, with the convention that this maximum is $0$ if no such metric space exists.

\item[(c)]
For any positive integers $p, q$, if $g(p-1, q) = g(p, q-1) = g(p-1, q-1) = 0$, then $g(p, q) = 0$. 

\end{itemize}
\end{lem}

\begin{proof}
(a)
The inequalities are clear from the definition and Lemma \ref{lem:finiteness}.

(b)
Such an $n$-point set $S = \{ v_0, \ldots, v_{n-1} \} \subset \R^{p, q}$ as in Lemma \ref{lem:finiteness} defines a metric $d$ on an $n$-point set $X = \{ 0, 1, \ldots, n-1 \}$ by $d(i, j) = \lVert v_i - v_j \rVert^2_{p, q}$ whose excess inertia is $(p, q, n - p - q - 1)$. Conversely, a metric $d$ on the $n$-point set $X$ whose excess inertia is $(p, q, n - p - q - 1)$ is realized as $d(i, j) = \lVert v_i - v_j \rVert^2_{p, q}$ by means of a non-obtuse set $S = \{ v_0, \ldots, v_{n-1} \} \subset \R^{p, q}$ such that $\Aff(S) = \R^{p, q}$ by Proposition \ref{prop:Schoenberg_embedding} and Proposition \ref{prop:metric_and_non_obtuse_set}. This correspondence allows the interpretation of $g(p, q)$.

(c)
Assume that $g(p, q) \neq 0$. Then, in view of (b), there exists an $n$-point metric space $(X, d)$ whose excess inertia is $(p, q, n - p - q - 1)$. Removing some points if necessary, we can assume that $n$ is the minimum. Thus, there exists no $(n-1)$-point metric subspace of $X$ whose excess inertia is $(p, q, n - p - q - 2)$. Now, let $Y$ be an $(n-1)$-point subspace of $X$. By the Cauchy interlace theorem, the excess inertia of $Y$ is either 
\begin{align*}
&(p - 1, q, n_0), &
&(p, q - 1, n_0), &
&(p - 1, q - 1, n_0 + 1), &
&(p, q, n_0 - 1),
\end{align*}
where $n_0 = n - p - q - 1$. By the assumption of (c) that $g(p - 1, q) = g(p, q - 1) = g(p - 1, q - 1) = 0$, the excess inertia of $Y$ is $(p, q, n - p - q - 2)$. However, this contradicts the present assumption of $X$. Therefore we conclude $g(p, q) = 0$.
\end{proof}

\begin{thm} \label{thm:maximal_points_in_non_obtuse_set_in_Rpq}
The following holds true.
\begin{itemize}
\item[(0)]
We have $g(0, 0) = 1$ and $g(0, q) = 0$ for $q \ge 1$.

\item[(1)]
We have $g(1, 0) = 2$ and $g(1, q) = 0$ for $q \ge 1$.

\item[(2)]
We have $g(2, 0) = 4$ and $g(2, q) = 0$ for $q \ge 1$.

\item[(3)]
We have $g(3, 0) = 8$, $g(3, 1) = 5$ and $g(3, q) = 0$ for $q \ge 2$.

\item[(4)]
We have $g(4, 0) = 16$, $g(4, 1) \ge 10$, $g(4, q) \ge 9$ for $q = 2, 3, 4$.

\end{itemize}
\end{thm}

The results in the theorem are summarized as follows:
$$
\begin{array}{|c|c|c|c|c|c|c|}
\hline
g(p, q) & q = 0 & q = 1 & q = 2 & q = 3 & q = 4 & q \ge 5 \\
\hline
p = 0 & 1 & 0 & 0 & 0 & 0 & 0 \\
\hline
p = 1 & 2 & 0 & 0 & 0 & 0 & 0 \\
\hline
p = 2 & 4 & 0 & 0 & 0 & 0 & 0 \\
\hline
p = 3 & 8 & 5 & 0 & 0 & 0 & 0 \\
\hline
p = 4 & 16 & \ge 10 & \ge 9 & \ge 9 & \ge 9 &  \\
\hline
\end{array}
$$
Note that $g(4, q)$ for $q \ge 5$ is undetermined. Note also that Theorem \ref{ithm:non_obtuse_problem} in \S\ref{sec:introduction} follows from the table above.

\begin{proof}
For $p \ge 0$, we have $g(p, 0) = 2^p$ by the theorem of Danzer-Gr\"{u}nbaum (Proposition \ref{prop:fact_non_obtuse_and_acute} (a)). Hence we only study the case where $q > 0$. 

(0)
For $q \ge 1$, we have $g(0, q) = 0$, because the inner product on $\R^{0, q}$ is negative definite. 

(1) 
Since there exists no $3$-point metric space of excess inertia $(1, 1, 0)$, we have $g(1, 1) = 0$. Then $g(1, q) = 0$ for $q \ge 2$ by Lemma \ref{lem:gpq_general_property} (c).

(2)
There exists no $4$-point metric space of excess inertia $(2, 1, 0)$. Hence $g(2, 1) = 0$. We then have $g(2, q) = 0$ for $q \ge 2$ by Lemma \ref{lem:gpq_general_property} (c).

(3)
There exists a $5$-point metric space of excess inertia $(3, 1, 0)$. Hence $g(3, 1) \ge 5$. There exists no $6$-point metric space of excess inertia $(3, 1, 1)$. It follows that $g(3, 1) = 5$. Because there exists no $6$-point metric space of excess inertia $(3, 2, 0)$, we have $g(3, 2) = 0$. Then $g(3, q) = 0$ for $q \ge 3$ by Lemma \ref{lem:gpq_general_property} (c).

(4) 
There exists a $10$-point metric space of excess inertia $(4, 1, 4)$, which implies $g(4, 1) \ge 10$. There are $9$-point metric spaces which have excess inertia $(4, 2, 2)$, $(4, 3, 1)$ and $(4, 4, 0)$. It follows that $g(4, 2) \ge 9$, $g(4, 3) \ge 9$ and $g(4, 4) \ge 9$.
\end{proof}

\subsection{The small-scale limit of magnitude}

We recall the definition of the magnitude (function) of a finite metric space \cite{L}.

\begin{dfn}
Let $X = (X, d)$ be an $n$-point metric space.
\begin{itemize}
\item
We define the \textit{zeta matrix} of $X$ as $Z_X = ( e^{-d_{x, y}} )_{x, y \in X} \in M(n, \R)$. 

\item
A \textit{weighting} is a vector $w = (w_x)_{x \in X} \in \R^n$ such that $Z_Xw = \1$, where $\1$ stands for the vector whose entries are $1$.

\item
If $X$ admits a weighting $w$, then the \textit{magnitude} of $X$ is defined as the sum of the entries of $w$
$$
\Mag(X) = \sum_{x \in X} w_x = \langle w, \1 \rangle.
$$

\item
For $t > 0$, we write $tX = (X, td)$ for $X$ equipped with the rescaled metric $td$ given by $(td)_{x, y} = t \cdot d_{x, y}$. The \textit{magnitude function} of $X$ is the function $t \mapsto \Mag(tX)$, whose domain consists of $t > 0$ such that $\Mag(tX)$ is defined. 

\end{itemize}
\end{dfn}

It is easy to see that the definition of the magnitude above is independent of the choice of a weighting. Note that $\Mag(tX)$ is defined for all sufficiently small $t > 0$ (see \cite{RY}, Remark 2.1). The limit $\lim_{t \to 0}\Mag(tX)$, if it exists, is called the \textit{small-scale limit} of the magnitude of $X$. We write $\SSL(X) = \lim_{t \to 0}\Mag(tX)$ for this limit.

As is mentioned in  \S\ref{sec:introduction}, we have $\SSL(X) = 1$ for any non-empty finite metric space $X$ consisting of at most four points. In the case of five points, the following result was first proved by Asao and Taira.

\begin{thm}[Asao-Taira] \label{thm:5_point_SSL}
We have $\SSL(X) \ge 1$ for any $5$-point metric space $X$ which admits a small-scale limit.
\end{thm}

The argument of the proof in the $5$-point case can be generalized to the $6$-point case, and we have the following result stated as Theorem \ref{ithm:SSL_6_point} in \S\ref{sec:introduction}.

\begin{thm} \label{thm:6_point_SSL}
We have $\SSL(X) \ge 1$ for any $6$-point metric space $X$ which admits a small-scale limit.
\end{thm}

\begin{proof}[Proof of Theorem \ref{thm:5_point_SSL} and Theorem \ref{thm:6_point_SSL}]
The following proof is based on an idea of Asao and Taira: According to the classifications of excess inertia of $5$- and $6$-point metric spaces in Theorem \ref{thm:5_point} and Theorem \ref{thm:6_point}, it suffices to see the following general results.
\begin{itemize}
\item[(a)]
We have $\SSL(X) \ge 1$ for any finite metric space $X$ whose excess inertia $(n_+, n_-, n_0)$ satisfies $n_- = 0$.

\item[(b)]
We have $\SSL(X) = 1$ for any finite metric space $X$ whose excess inertia $(n_+, n_-, n_0)$ satisfies $n_0 = 0$.

\end{itemize}

(a)
As is mentioned in Remark \ref{rem:negative_type}, the condition $n_- = 0$ means that $X$ is of negative type. As shown in \cite{M}, if $X$ is of negative type, then $Z_{tX}$ is positive definite for all $t > 0$. It is known \cite{L} that the positive definiteness of $Z_{tX}$ implies $\Mag(tX) \ge 1$ for all $t > 0$. Hence $\SSL(X) = \lim_{t \to 0}\Mag(tX) \ge 1$.

(b)
In \cite{RY}, a quantity $F_X \in \R$ is introduced for any finite metric space $X$, which has the property that we have $\SSL(X) = 1$ whenever $F_X \neq 0$. An expression of $F_X$ is $F_X = (-1)^{n-1} \langle \1, \mathrm{adj}(D) \1 \rangle$, where $\mathrm{adj}(D)$ is the adjugate of the distance matrix $D \in M(n, \R)$ of $X$, which is assumed to consist of $n$ points. Let us define  $\mathrm{CM} \in M(n+1, \R)$ by 
$$
\mathrm{CM} 
=
\left(
\begin{array}{cc}
0 & {}^t\1 \\
\1 & D
\end{array}
\right).
$$
Then $\det \mathrm{CM} = -\langle \1, \mathrm{adj}(D) \1 \rangle $ by a cofactor expansion. Consider the following expression of $\mathrm{CM}$ and an invertible matrix $P \in M(n+1, \R)$
\begin{align*}
\mathrm{CM} 
&=
\left(
\begin{array}{cc}
0 & {}^t\1 \\
\1 & D
\end{array}
\right)
=
\left(
\begin{array}{ccc}
0 & 1 & {}^t \1 \\
1 & 0 & {}^t d^{(0)} \\
\1 & d^{(0)} & D^{(0)} 
\end{array}
\right),
&
P
&=
\left(
\begin{array}{ccc}
1 & 0 & - {}^td^{(0)} \\
0 & 1 & - {}^t\1 \\
0 & 0 & I
\end{array}
\right),
\end{align*}
where $I \in M(n-1, \R)$ is the identity matrix. By means of $d^{(0)} \in \R^{n-1}$ and $D^{(0)} \in M(n-1, \R)$, we can express the $0$th excess matrix as $E^{(0)} = \1{}^td^{(0)} + d^{(0)}{}^t\1 - D^{(0)}$. We can then verify that
$$
{}^t P (\mathrm{CM}) P
=
\left(
\begin{array}{ccc}
0 & 1 & 0 \\
1 & 0 & 0 \\
0 & 0 & -E^{(0)}
\end{array}
\right).
$$
It follows that $\det \mathrm{CM} = (-1)^n \det E^{(0)}$. Now we have
$$
F_X 
= (-1)^{n-1} \langle \1, \mathrm{adj}(D) \1 \rangle
= (-1)^n \det \mathrm{CM}
= \det E^{(0)}.
$$
Thus, if $X$ has the excess inertia $(n_+, n_-, 0)$, then $\SSL(X) = 1$.
\end{proof}


\begingroup
\raggedright

\endgroup

\end{document}